\documentclass[11pt]{article}

\usepackage[margin=1in]{geometry}
\usepackage[utf8]{inputenc}
\usepackage[T1]{fontenc}
\usepackage{amssymb,latexsym,amsmath,amsthm,mathtools,graphicx,enumerate}
\usepackage{tabularx,booktabs,longtable,xcolor,siunitx,url}
\usepackage[nameinlink]{cleveref}
\allowdisplaybreaks
\providecommand{\texorpdfstring}[2]{#1}

\newtheorem{theorem}{Theorem}
\newtheorem{lemma}{Lemma}
\newtheorem{proposition}{Proposition}
\newtheorem{corollary}{Corollary}
\theoremstyle{definition}
\newtheorem{definition}{Definition}
\newtheorem{conjecture}{Conjecture}
\newtheorem{remark}{Remark}
\newtheorem{example}{Example}

\newcommand{\height}{H}

\newcommand{\p}{\boldsymbol p}

\title{Completely Additive Height Functions: Profile Laws, Matula Bounds, and Inverse Growth}
\author{Hartosh Singh Bal\\
\small The Caravan, New Delhi, India\\
\small\texttt{hartoshbal@gmail.com}}
\date{}

\begin{document}
\maketitle

\begin{abstract}
The height $H(n)$ of an integer $n$ is classically the number of iterations of
Euler's totient function required to reach $1$. H.~N.~Shapiro showed that a
modification of this function is completely additive.  We study completely
additive height functions with finite prime fibers.  Their prime-height profile
$\pi_k$ determines the height multiplicities $N_k$ through the
weighted-multipartition identity
$\sum_k N_k q^k=\prod_j(1-q^j)^{-\pi_j}$, and conversely every profile
containing infinitely many primes is realizable.  We introduce iteratively
defined heights encompassing Shapiro-type totient heights and the Matula height.
For the Matula height, we give purely number-theoretic proofs of the classical Gutman–Ivić extremal bounds, thereby answering their question whether the maximal bound can be derived without recourse to the rooted-tree interpretation.  Using Meinardus' theorem in its full form, we prove a
conditional inverse-growth law: if $\Pi_k\sim Ck^\alpha$, then
$\log N_k\sim C_2 k^{\alpha/(\alpha+1)}$, with an explicit constant.  We also
derive average-order results for a canonical sequential realization and report
computations for the Shapiro height beyond the polynomial regime.
\end{abstract}

\section{Introduction}

A \emph{height function} in this paper is a completely additive function
$H:\mathbb{N}\to\mathbb{N}_0$ with $H(1)=0$ and $H(p) \ge 1$, such that for each $n\ge1$ there
are only finitely many primes $p$ with $H(p)=n$.
Writing
\[
  \pi_n := \#\{p\text{ prime}: H(p)=n\}
  \qquad\text{and}\qquad
  N_n := \#\{m\ge1:H(m)=n\},
\]
complete additivity immediately forces the weighted-partition identity
\[
  \sum_{n\ge0} N_n q^n
  \;=\;
  \prod_{k\ge1} (1-q^k)^{-\pi_k}.
\]
Thus, the prime-height profile $(\pi_k)$ canonically determines the height
multiplicities $(N_n)$ via a multipartition Euler product with weights
$b_k=\pi_k$.  Conversely, given any sequence $(\pi_k)_{k\ge1}$ of nonnegative
integers with $\sum_{k\ge1}\pi_k=\infty$, one can realize it as the prime-height
profile of some height function (by assigning heights to the primes accordingly
and extending by complete additivity; see Proposition~\ref{prop:realize-profile}).
In this setting the study of height functions with finite prime fibers aligns
naturally with multipartitions and their asymptotics.  The benchmark profiles
$\pi_k\equiv 1$ and $\pi_k=k$ recover ordinary and plane partitions.

The paper develops this equivalence in two complementary directions.  On the
structural side, simple recursive prescriptions for $H(p_i)$ produce rich
combinatorial classes (most notably rooted trees via the Matula height) and also
capture heights arising from iterates of multiplicative functions.  On the
analytic side, the prime-height profile controls both the average order of
$H(n)$ and the subexponential growth of the induced multipartition layers,
placing natural examples into a small number of growth classes.

\subsection{Main Results}

The weighted-partition identity used below belongs to a classical circle of
ideas around Euler products, additive functions, and multipartitions.  The
contribution of the present paper is to make this identity canonical for
completely additive height functions with finite prime fibers, and then to use it
as an organizing principle for recursive prime-height constructions, Matula
bounds, and conditional inverse-growth laws.

The principal results fall into three groups.

First, for a completely additive height function $H$, the prime-height profile
$\pi_k=\#\{p:H(p)=k\}$ determines the height multiplicities by
\[
  \sum_{n\ge0}N_nq^n=\prod_{k\ge1}(1-q^k)^{-\pi_k}.
\]
We also give a general recursive construction theorem of the form
$H(p_i)=H(f(i))+1$, encompassing the Matula and Shapiro-type examples.

Second, for the Matula height, we give number-theoretic proofs of the classical
Gutman--Ivi\'c extremal bounds for the smallest and largest integers at fixed
height, using only complete additivity, the recursion $H(p_i)=H(i)+1$, and
explicit estimates for the $i$th prime.

Third, for profiles with polynomial cumulative growth, the profile supplies the
pole and residue of the associated Dirichlet series.  Under the full hypotheses
of Meinardus' theorem, including the minor-arc condition, this yields a precise
stretched-exponential law for $N_n$.  The statement is conditional precisely
because a power law for $\Pi(x)$ alone does not rule out lattice support.  We
also study the canonical sequential realization of a prescribed profile and its
average order.

Our first structural theme is the class of \emph{iteratively defined} heights.
In Section~\ref{sec3} we consider heights determined on primes by
\[
  H(p_i) = H(f(i)) + 1,
\]
where $p_i$ is the $i$th prime and $f:\mathbb{N}\to\mathbb{N}$ satisfies a
natural finiteness and ``no-backtracking'' condition.  Theorem~\ref{thm:iterative}
shows that this covers a large class of examples, including the Matula height
$H(p_i)=H(i)+1$ and Shapiro-type heights arising from iterates of Euler's
$\varphi$-function.  In the Matula case we revisit a result of Gutman and Ivi\'c~\cite{gutman} on
extremal Matula numbers at fixed height.  They observed that their proof of the
upper bound for the maximal Matula number is graph-theoretic and remarked that
``it is not clear'' how to deduce it by purely number-theoretic arguments.  In
\Cref{sec3} we supply such an arithmetic derivation for \emph{both} extremal
bounds (minimal and maximal), using only complete additivity, the recursion
$H(p_i)=H(i)+1$, and explicit prime-index bounds, with no appeal to the rooted-tree
interpretation.

Our main analytic theme is an \emph{inverse-growth} principle in the polynomial,
profile-driven regime.  Writing
\[
  B_H(s) := \sum_{k\ge1} \frac{\pi_k}{k^s},
\]
Section~\ref{sec6} gives the analytic inverse-growth theorem.  The section is
formulated conditionally under the full hypotheses of Meinardus' theorem: the
profile growth gives the dominant pole and residue, while the minor-arc
condition is an additional hypothesis and must be checked separately.  In this
setting one obtains
\[
  \log N_n \sim C_2n^{\alpha/(\alpha+1)}
  \qquad (n\to\infty),
\]
with
\[
  C_2=
  \left(1+\frac1\alpha\right)
  \bigl(C\Gamma(\alpha+1)\zeta(\alpha+1)\bigr)^{1/(\alpha+1)}
\]
when the residue at the pole is $C$.  In Section~\ref{sec7} we develop this analysis
further: given a polynomial profile we define a canonical \emph{sequential
realization} $\height_{\mathrm{seq}}$ and derive corresponding average-order
laws for $\sum_{n\le x}\height_{\mathrm{seq}}(n)$ using the additive-function
reduction proved in Section~\ref{sec5}.

A secondary motif is the squarefree, or radical-class, variant of the theory.
Here repeated prime factors are suppressed from the outset, so that the Euler
product changes from an unrestricted multipartition product to a distinct-part
product.  For non-recursive prime-height assignments this agrees with the
literal radicalization \(n\mapsto \operatorname{rad}(n)\), but for recursive
height structures we impose the recursion directly on squarefree representatives.
This produces distinct-part analogs of the preceding constructions, including
Shapiro--Erd\H{o}s and Dedekind-type examples.  We develop
this in Section~\ref{sec4}.

Finally, we include a brief experimental discussion of the Shapiro totient
height to illustrate that recursive heights can impose strong arithmetic
structure on the prime profile and may leave the polynomial regime.  This is
meant as evidence that substantially richer behavior can occur beyond the
universal polynomial picture.

The paper is organized as follows.
Section~\ref{sec2} recalls the basic notion of a height function and the induced
multipartition identity.
Section~\ref{sec3} develops the iteratively defined heights $H(p_i)=H(f(i))+1$
and proves the Matula extremal bounds in a purely number-theoretic way.
Section~\ref{sec4} treats squarefree height functions and the resulting
distinct-part structures.
Section~\ref{sec5} discusses average order and proves the additive-function
reduction used later.
Section~\ref{sec6} proves the inverse-growth theorem via Meinardus' method, and
Section~\ref{sec7} derives average-order consequences in the polynomial profile
regime for the canonical sequential realization.
Section~\ref{sec8} presents computational observations for the Shapiro height.
We conclude with open problems and conjectural directions arising from the
height perspective.

\section{A Class of Completely Additive Functions}\label{sec2}

In this section we record some general facts about completely additive
functions and height structures, and illustrate them with examples which realize
familiar partition families.  Throughout, a \emph{height function} $H$ is a
completely additive function
\[
  H:\mathbb{N}\to\mathbb{N}_0
\]
such that $H(1)=0$, such that $H(p)\ge1$ for every prime $p$, and such that, for
each $n\ge1$, the set of primes
\[
  \{p : H(p)=n\}
\]
is finite.  We write $\pi_n$ for the number of primes at height $n$.

It will be convenient to formalize the combinatorial object encoded by
a height function.

\begin{definition}\label{def1}
Given a sequence $(\pi_n)_{n\ge1}$ of nonnegative integers, a
\emph{multipartition} of $n$ with $\pi_\ell$ colors at part size
$\ell$ is a family of nonnegative integers
$\alpha_{\ell,j}$ indexed by $1\le\ell\le n$ and
$1\le j\le \pi_\ell$ such that
\[
  \sum_{\ell=1}^n \ell \sum_{j=1}^{\pi_\ell} \alpha_{\ell,j}
  \;=\; n.
\]
The generating function for the numbers $M_n$ of such multipartitions
is
\[
  \sum_{n\ge0} M_n q^n
  \;=\;
  \prod_{\ell\ge1} \frac{1}{(1-q^\ell)^{\pi_\ell}}.
\]
\end{definition}

\begin{proposition}\label{thm:height-multipartition}
A height function $H$ induces a bijection between $\mathbb{N}$ and
the disjoint union of all multipartitions, sending $1$ to the empty multipartition, with the multipartition parameters $\pi_n$ given by the prime heights of $H$.
\end{proposition}

\begin{proof}
By hypothesis, for each $n\ge1$ the number
$\pi_n := \#\{p : H(p)=n\}$ is finite.  Enumerate the primes at height
$n$ as $p_{n,1},\dots,p_{n,\pi_n}$.  Let $m\in\mathbb{N}$, and write
its prime factorization as
\[
  m
  \;=\;
  \prod_{n\ge1}\;\prod_{j=1}^{\pi_n} p_{n,j}^{\alpha_{n,j}},
\]
where only finitely many of the exponents $\alpha_{n,j}$ are nonzero.
By complete additivity,
\[
  H(m)
  \;=\;
  \sum_{n\ge1}\;\sum_{j=1}^{\pi_n} \alpha_{n,j} H(p_{n,j})
  \;=\;
  \sum_{n\ge1} n \sum_{j=1}^{\pi_n} \alpha_{n,j}.
\]
Thus, for each fixed $k\ge0$, the condition $H(m)=k$ is equivalent to
\[
  \sum_{n=1}^k n \sum_{j=1}^{\pi_n} \alpha_{n,j} = k.
\]
In other words, the exponent data
$\bigl(\alpha_{n,j}\bigr)_{1\le n\le k, 1\le j\le\pi_n}$ is exactly a
multipartition of $k$ with $\pi_n$ generators at each size $n$.
Conversely, given such data, unique factorization into primes produces
a unique integer $m$ with $H(m)=k$.  Taking the disjoint union over
all $k\ge0$ yields a bijection between $\mathbb{N}$ (via $H$) and
the set of all multipartitions.
\end{proof}

As a consequence, the level-size sequence $N_n := \#\{m:H(m)=n\}$ is
given by the generating function
\begin{equation}\label{eq:GN}
  \sum_{n\ge0} N_n q^n
  \;=\;
  \prod_{n\ge1} \frac{1}{(1-q^n)^{\pi_n}}.
\end{equation}

\begin{proposition}[Realizing a prescribed profile]\label{prop:realize-profile}
Let $(\pi_k)_{k\ge1}$ be a sequence of nonnegative integers with $\sum_{k\ge1}\pi_k=\infty$.
Then there exists a height function $H$ such that
$\#\{p\text{ prime}:\, H(p)=k\}=\pi_k$ for every $k\ge1$.
\end{proposition}

\begin{proof}
Enumerate the primes as $p_1<p_2<\cdots$, and set
\[
  S_k=\sum_{\ell=1}^k\pi_\ell,\qquad S_0=0.
\]
Since $\sum_k\pi_k=\infty$, the intervals
\[
  I_k=\{i\ge1:S_{k-1}<i\le S_k\}
\]
partition $\mathbb N$, with $|I_k|=\pi_k$.  Define $H(p_i)=k$ for
$i\in I_k$, set $H(1)=0$, and extend to all of $\mathbb N$ by complete
additivity.  Then $H(p)\ge1$ for every prime $p$, and by construction exactly
$\pi_k$ primes have height $k$.  Since each $\pi_k$ is finite, $H$ is a height
function.
\end{proof}

We next record a simple sufficient condition ensuring that a completely
additive function is a height function.

\begin{proposition}\label{prop:f-static}
Let $f:\mathbb{N}\to\mathbb{N}$ be an arithmetic function such that
for each $t\in\mathbb{N}$, the equation $f(i)=t$ has only finitely
many solutions $i\in\mathbb{N}$.  Define a completely additive
function $H$ by
\[
  H(1)=0,\qquad H(p_i)=f(i),\qquad
  H\Bigl(\prod p_i^{\alpha_i}\Bigr)=\sum_i \alpha_i H(p_i).
\]
Then $H$ is a height function, with $\pi_n = \#\{i : f(i)=n\}$.
\end{proposition}

\begin{proof}
By construction $H$ is completely additive and $H(1)=0$.  For each
$n\ge1$,
\[
  \pi_n =
  \#\{p : H(p)=n\}
  = \#\{i : f(i)=n\}<\infty
\]
by hypothesis on $f$.  Thus, $H$ satisfies the defining conditions of a
height function.
\end{proof}

We also record an algorithmic description of the level sets
$H^{-1}(n)$, assuming that the primes at each lower level are already
known.

\begin{theorem}[Level construction algorithm]\label{thm:algo}
Suppose that for some $n\ge1$ the sets
\[
  \mathcal{H}_k := \{m\in\mathbb{N} : H(m)=k\}
\]
and the subsets of primes $\mathcal{P}_k:=\{p\in\mathcal{H}_k\}$ are
known for all $1\le k\le n-1$.  Then the set
$\mathcal{H}_n=\{m:H(m)=n\}$ is obtained as follows:
\begin{enumerate}
\item For each $\ell$ with $1\le \ell \le \lfloor n/2\rfloor$, form
  the products $p\cdot m$ with $p\in\mathcal{P}_\ell$ and
  $m\in\mathcal{H}_{n-\ell}$.
\item Let $\mathcal{C}_n$ be the set of all such products, with
  duplicates removed.
\item Let $\mathcal{P}_n$ be the set of primes at height $n$.  Then
  \[
    \mathcal{H}_n = \mathcal{C}_n \cup \mathcal{P}_n.
  \]
\end{enumerate}
\end{theorem}

\begin{proof}
Let $m\in\mathcal{H}_n$ be composite, and write its prime factorization as
\[
  m=\prod_{i=1}^t p_i^{\alpha_i}.
\]
Since $m$ is composite, the total number of prime factors counted with
multiplicity is
\[
  \Omega(m):=\sum_{i=1}^t\alpha_i\ge2.
\]
Complete additivity gives
\[
  n=H(m)=\sum_{i=1}^t \alpha_iH(p_i),
\]
with every summand $H(p_i)$ a positive integer.  Hence, at least one prime
divisor $p$ of $m$ satisfies $H(p)\le n/2$; otherwise every prime factor counted
with multiplicity would contribute more than $n/2$, and the total would exceed
$n$.  Choose such a prime $p$ and put
$\ell:=H(p)\le \lfloor n/2\rfloor$.  Then
\[
  H(m/p)=H(m)-H(p)=n-\ell,
\]
so $m/p\in\mathcal{H}_{n-\ell}$ and $m$ is produced in Step~(1).

Conversely, if $p\in\mathcal{P}_\ell$ and $m'\in\mathcal{H}_{n-\ell}$,
then complete additivity gives
\[
  H(pm') = H(p)+H(m') = \ell + (n-\ell) = n,
\]
so every product produced in Step~(1) lies in $\mathcal{H}_n$.  Finally
the primes of height $n$ are precisely $\mathcal{P}_n$ by definition,
so $\mathcal{H}_n$ is the disjoint union of $\mathcal{C}_n$ (composite
elements) and $\mathcal{P}_n$ (prime elements).
\end{proof}

We now turn to concrete choices of $H$.

\subsection{Ordinary Partitions}
\label{ex:H=i}

Let $H(p_i)=i$ for each $i\ge1$.  Then
\[
  H\Bigl(\prod_{i=1}^r p_i^{\alpha_i}\Bigr)
  \;=\; \sum_{i=1}^r i\,\alpha_i
\]
and there is exactly one prime at each height $n$, namely $p_n$.
Hence, $\pi_n=1$ for all $n\ge1$, so by Equation~\eqref{eq:GN}
\[
  \sum_{n\ge0} N_n q^n
  \;=\;
  \prod_{n\ge1} \frac{1}{1-q^n},
\]
and $N_n$ is the ordinary partition function $p(n)$.

Table~\ref{tab:H=i} shows the first few levels.

\begin{table}[ht]
\centering
\begin{tabular}{| l | l | }
\hline\hline
$n$ & $\{m\mid H(m)=n\}$\\
\hline
$6$ & $\mathbf{13}, 21, 22, 25, 27, 28, 30, 36, 40, 48, 64$ \\
$5$ & $\mathbf{11}, 14, 15, 18, 20, 24, 32$ \\
$4$ & $\mathbf{7}, 9, 10, 12, 16$ \\
$3$ & $\mathbf{5}, 6, 8$ \\
$2$ & $\mathbf{3}, 4$ \\
$1$ & $\mathbf{2}$ \\
$0$ & $1$ \\
\hline
\end{tabular}
\caption{Numbers with height $\le 6$ for $H(p_i)=i$. Primes are in bold.}
\label{tab:H=i}
\end{table} 

The combinatorial correspondence is completely transparent.  For
instance, the $7$ elements at height $5$,
\[
  11, 14, 15, 18, 20, 24, 32,
\]
correspond bijectively to the $7$ partitions of $5$:
\begin{itemize}
\item $11$ is prime, so corresponds to the partition $5$;
\item $14=2\cdot7$ corresponds to $1+4$;
\item $15=3\cdot5$ corresponds to $2+3$;
\item $18=2\cdot 3^2$ corresponds to $1+2+2$;
\item $20=2^2\cdot 5$ corresponds to $1+1+3$;
\item $24=2^3\cdot 3$ corresponds to $1+1+1+2$;
\item $32=2^5$ corresponds to $1+1+1+1+1$.
\end{itemize}

The largest element at height $n$ is easy to identify.

\begin{proposition}\label{prop:max-H=i}
For the height function $H(p_i)=i$, the largest number at height $n$
is $2^n$.
\end{proposition}

\begin{proof}
Let $m$ have height $n$ and prime factorization
$m=\prod_{i\ge1}p_i^{\alpha_i}$.  Then
\[
  n=H(m)=\sum_{i\ge1}i\alpha_i.
\]
Since $p_i\le 2^i$ for every $i\ge1$, we have
\[
  m=\prod_{i\ge1}p_i^{\alpha_i}
  \le \prod_{i\ge1}2^{i\alpha_i}
  =2^{\sum_i i\alpha_i}
  =2^n.
\]
Equality can occur only if $p_i=2^i$ whenever $\alpha_i>0$.  Since this happens
only for $i=1$, equality forces $\alpha_1=n$ and $\alpha_i=0$ for all $i\ge2$.
Thus, $m=2^n$, and it is the unique maximal element at height $n$.
\end{proof}

Identifying the smallest element at height $n$ is more delicate.
Numerically it appears that $p_n$ is always the smallest element at
height $n$, and this can indeed be proved using explicit bounds for
the primes and the inequality $p_a p_b \ge p_{a+b}$ for all
$a,b\ge1$. The proof replicates the argument carried out later in the paper in the case of the Matula function \cite{gutman}.

\begin{remark}
Using explicit bounds for $p_n$ due to Rosser--Schoenfeld, one can
show that for all $a,b\ge1$,
\[
  p_a p_b \;\ge\; p_{a+b}.
\]
It follows that if $H(p_i)=i$ and $H(m)=n$, then any factorization
$n=i_1+\dots+i_k$ leads to
\[
  \prod_{r=1}^k p_{i_r} \;\ge\; p_{i_1+\dots+i_k}=p_n,
\]
with equality only for the trivial partition $n=n$.  Thus, $p_n$ is the
unique smallest element at height $n$.  We will not need this fact in
what follows.
\end{remark}

\subsection{Prime Partitions}
\label{ex:H=prime}

Let $H(p_i)=p_i$ for each $i\ge1$.  Then
\[
  H\Bigl(\prod_{i=1}^r p_i^{\alpha_i}\Bigr)
  \;=\; \sum_{i=1}^r \alpha_i p_i.
\]
In this case there is exactly one prime at each \emph{prime} height,
and no primes at composite heights.  More precisely,
\[
  \pi_n = \begin{cases} 1, & n\ \text{prime},\\[3pt]
                        0, & \text{otherwise}.
          \end{cases}
\]
Therefore, Equation~\eqref{eq:GN} becomes
\[
  \sum_{n\ge0} N_n q^n
  \;=\;
  \prod_{n\ge1} \frac{1}{(1-q^n)^{\chi(n)}},
\]
where $\chi(n)$ is the characteristic function of the primes.  Thus,
$N_n$ is the number of partitions of $n$ into prime parts.  Table~\ref{tab:H=prime} shows the first few levels.

\begin{table}[ht]
\centering
\begin{tabular}{| l | l | }
\hline\hline
$n$ & $\{m\mid H(m)=n\}$\\
\hline
$8$ & $15, 16, 18$ \\
$7$ & $\mathbf{7}, 10, 12$ \\
$6$ & $8, 9$ \\
$5$ & $\mathbf{5}, 6$ \\
$4$ & $4$ \\
$3$ & $\mathbf{3}$ \\
$2$ & $\mathbf{2}$ \\
$1$ & \\
$0$ & $1$ \\
\hline
\end{tabular}
\caption{Numbers with height $\le 8$ for $H(p_i)=p_i$. Primes are in bold.}
\label{tab:H=prime}
\end{table} 

Here again one can study extremal elements.  The smallest element at
height $n$ is at least $p_{\pi(n)}$, while arguments about
maximizing products with fixed sum of exponents show that for $n\ge2$
the largest element is given by
\[
  3^{n/3},\quad 4\cdot 3^{(n-4)/3},\quad
  2\cdot 3^{(n-2)/3}
\]
for $n\equiv0,1,2\pmod{3}$, respectively.  We omit the
details, as this is worked out in the paper by Alladi and Erd\H{o}s \cite{alladi}.

\subsection{Plane Partitions}
\label{ex:plane}

A particularly interesting special case of
Proposition~\ref{prop:realize-profile} is obtained by taking
$\pi_n=n$ for all $n$.  For concreteness, we may define
\[
  H(p_i)=n
  \quad\text{whenever}\quad
  \frac{n(n-1)}{2} < i \le \frac{n(n+1)}{2},
\]
so that there are exactly $n$ primes at height $n$.  Then
\[
  \sum_{n\ge0} N_n q^n
  \;=\;
  \prod_{n\ge1} \frac{1}{(1-q^n)^n}.
\]
The product on the right is the classical generating function for
plane partitions of $n$.  Thus, \(N_n\) is both the number of multipartitions of \(n\) with \(k\) colors at part size \(k\), and the number of plane partitions
of \(n\); see, for instance, Andrews~\cite{andrews}.

Although the level sizes coincide with those of plane partitions, the
map $m\mapsto H(m)$ does not directly preserve the natural additive
structure on plane partitions, and the identification a priori depends
on knowing that the two counting sequences coincide.  In later
sections we will see examples where the height structure gives a more
direct and constructive bridge between integers and higher-dimensional
partition objects.

Table~\ref{tab:plane} shows the first few levels.

\begin{table}[ht]
\centering
\begin{tabular}{| l | l | }
\hline\hline
$n$ & $\{m\mid H(m)=n\}$\\
\hline
$5$ & $18, 21, 24, 28, 30, \mathbf{31}, 32, 33, 34, 35, \mathbf{37}, 38, 39, 40,$\\
    & $\mathbf{41}, \mathbf{43}, 44, 46, \mathbf{47}, 50, 52, 55, 58, 65 $\\
$4$ & $9, 12, 14, 15, 16, \mathbf{17}, \mathbf{19}, 20, 22, \mathbf{23}, 25, 26, \mathbf{29}$ \\
$3$ & $6, \mathbf{7}, 8, 10, \mathbf{11}, \mathbf{13}$ \\
$2$ & $\mathbf{3}, 4, \mathbf{5}$ \\
$1$ & $\mathbf{2}$ \\
$0$ & $1$ \\
\hline
\end{tabular}
\caption{Numbers with height $\le 5$ for $\pi_n=n$. Primes are in bold.}
\label{tab:plane}
\end{table}

We end this section by noting that even simple-looking choices of $H$
can give rise to multipartition structures that are not easily
recognizable as classical partition families without explicit
computation of the prime heights.  One such example is
$H(p_i)=\lfloor p_i/i\rfloor$, for which the induced weights
$(\pi_n)_{n\ge1}$ grow in an irregular fashion and do not appear to
match any standard partition model.  We will not pursue this example
further here, but it illustrates that the height framework naturally
produces families well beyond the classical catalogs.

\section{Iteratively Defined Completely Additive Functions}
\label{sec3}

The description of $H$ in Section~\ref{sec2} suffices when the height
of each prime is prescribed in advance, so that the sequence
$(\pi_n)$ is known.  Many of the most interesting examples, however,
are \emph{iteratively defined}, with $H(p_i)$ depending on the height
of integers smaller than $p_i$.  In this section we give a general
framework for such definitions and show that they still produce height
functions, before specializing to concrete cases such as the Matula
numbers and Shapiro's totient height.

\begin{theorem}\label{thm:iterative}
Let $f:\mathbb{N}\to\mathbb{N}$ be an arithmetic function with
$f(1)=1$ such that:
\begin{enumerate}
\item for any prime $p_r$, if $p_r\mid f(i)$ then $r<i$;
\item for each $t\in\mathbb{N}$, the equation $f(i)=t$ has only
  finitely many solutions $i\in\mathbb{N}$.
\end{enumerate}
Fix an integer $j\ge1$ and define $H$ on primes by
\[
  H(2)=j,\qquad H(p_i) = H(f(i)) + j \quad (i\ge1),
\]
and extend $H$ to $\mathbb{N}$ by complete additivity and $H(1)=0$.
Then $H$ is a height function.
\end{theorem}
\begin{proof}
We treat the case $j=1$; the general case is identical (all heights are
shifted by the factor $j$).  We first check that the recursion defining
$H$ on primes is well-defined.

Define $H(2)=1$.  Suppose inductively that $H(p_r)$ has been defined for
all $r<i$.  To define $H(p_i)$ we must evaluate $H(f(i))$.
By hypothesis \textup{(1)}, every prime divisor of $f(i)$ is of the form
$p_r$ with $r<i$.  Hence, $H$ is already defined on all prime factors of
$f(i)$, and by complete additivity $H(f(i))$ is determined.  We may
therefore set
\[
H(p_i):=H(f(i))+1,
\]
and this determines $H$ on all primes by induction on $i$.

We next prove that the prime layers are finite.  Extend $H$ to $\mathbb N$ by complete additivity and $H(1)=0$.
We prove by induction on $n\ge1$ that the set of primes
\[
P_n:=\{p:\ H(p)=n\}
\]
is finite.

For $n=1$, we have $H(p_i)=1$ if and only if $H(f(i))=0$, i.e.\ $f(i)=1$.
By hypothesis \textup{(2)}, the equation $f(i)=1$ has only finitely many
solutions $i$, hence $P_1$ is finite.

Assume $P_\ell$ is finite for all $\ell<n$.  Let
\[
S_{n-1}:=\bigcup_{\ell=1}^{n-1} P_\ell
\]
be the finite set of primes of height at most $n-1$.
If $m\in\mathbb N$ satisfies $H(m)=n-1$, then in its prime factorization
$m=\prod_{p}p^{v_p(m)}$ only primes $p\in S_{n-1}$ can occur: any prime
with height $\ge n$ would contribute at least $n$ to $H(m)$.
Moreover, for each $p\in S_{n-1}$ we have $v_p(m)\le n-1$ since
$H(p)\ge 1$ and $v_p(m)H(p)\le H(m)=n-1$.
Therefore, the set
\[
\mathcal H_{n-1}:=\{m:\ H(m)=n-1\}
\]
is finite (it is contained in the finite set of products
$\prod_{p\in S_{n-1}} p^{e_p}$ with $0\le e_p\le n-1$).

Now if $H(p_i)=n$, then by definition $H(f(i))=n-1$, hence
$f(i)\in\mathcal H_{n-1}$.
For each fixed $a\in\mathcal H_{n-1}$, hypothesis \textup{(2)} implies
that the equation $f(i)=a$ has only finitely many solutions $i$.
Since $\mathcal H_{n-1}$ is finite, only finitely many indices $i$
satisfy $f(i)\in\mathcal H_{n-1}$, and thus $P_n$ is finite.

This completes the induction, so every prime-height layer $P_n$ is
finite.  As shown above, finiteness of the prime layers up to height $n$
implies finiteness of $\mathcal H_n$, hence each integer height layer is
finite as well.  Therefore, $H$ is a height function.
\end{proof}

This iterative framework encompasses a number of classical constructions.

\subsection{Matula Numbers}\label{ex:Matula}

Let $f(i)=i$ and $j=1$ in Theorem~\ref{thm:iterative}.  Then
\[
  H(p_i) = H(i) + 1.
\]
In particular,
\[
  H(2)=1,\quad H(3)=2,\quad H(5)=3,\quad H(7)=3,\quad H(11)=4,\quad\dots.
\]
The primes at each height $n$ are in bijection with the integers at
height $n-1$, so if we denote by $N_n$ the number of integers at
height $n$ and by $\pi_n$ the number of primes at height $n$, we have
\[
  \pi_n = N_{n-1}\qquad (n\ge1).
\]
Theorem~\ref{thm:height-multipartition} and Equation~\eqref{eq:GN} then yield
the recursive generating function
\begin{equation*}
  \sum_{n\ge0} N_n q^n
  \;=\;
  \prod_{n\ge1} \frac{1}{(1-q^n)^{N_{n-1}}}.
\end{equation*}

These $N_n$ count unlabeled rooted trees on $n+1$ vertices;
see, for example, Otter~\cite{otter}, or Gutman and Ivi\'c~\cite{gutman}.
Equivalently, the integers $m$ in $\mathcal{H}_n$ are the \emph{Matula numbers}
of rooted trees with $n+1$ vertices.

The multiplicative structure of $\mathcal{H}_n$ encodes the rooted tree
corresponding to a Matula number via the classical Matula map.
Let $\mathcal{T}$ be the set of unlabeled rooted trees and write $M:\mathcal{T}\to\mathbb N$
for the Matula number.  It is defined recursively by
\[
  M(\text{single vertex})=1,
  \qquad
  M(T)=\prod_{u\ \text{child of the root of }T} p_{\,M(T_u)},
\]
where $T_u$ denotes the rooted subtree hanging from the child $u$.
Equivalently, if
\[
  m=\prod_{r} p_r^{e_r},
\]
then $m$ encodes a rooted tree whose root has exactly $e_r$ branches of
type $T_r$, where $T_r$ is the rooted tree with Matula number $r$; iterating
this rule recovers the entire tree.

Table~\ref{tab:Matula} lists the first few levels.

\begin{table}[ht]
\centering
\begin{tabular}{| l | l | }
\hline\hline
$n$ & $\{m\mid H(m)=n\}$\\
\hline
$5$ & $15, 18, 20, 21, 22, \mathbf{23}, 24, 26, 28, \mathbf{29}, \mathbf{31}, 32, 34, \mathbf{37}, 38,$\\
    & $\mathbf{41}, \mathbf{43}, \mathbf{53}, \mathbf{59}, \mathbf{67}$ \\
$4$ & $9, 10, \mathbf{11}, 12, \mathbf{13}, 14, 16, \mathbf{17}, \mathbf{19}$ \\
$3$ & $\mathbf{5}, 6, \mathbf{7}, 8$ \\
$2$ & $\mathbf{3}, 4$ \\
$1$ & $\mathbf{2}$ \\
$0$ & $1$ \\
\hline
\end{tabular}
\caption{Numbers with height $\le 5$ in the Matula structure. Primes are in bold.}
\label{tab:Matula}
\end{table}

Gutman and Ivi\'c~\cite{gutman} studied extremal Matula numbers at
fixed height $n$, obtaining sharp asymptotic bounds for the smallest
and largest numbers in $\mathcal{H}_n$.  Their proofs use structural
properties of rooted trees.  We now show that the same extremal
results admit a completely additive proof, using only the recursive
definition $H(p_i)=H(i)+1$ and explicit bounds for $p_n$.

\subsection{Extremal Matula Numbers after Gutman--Ivi\'c}\label{gutman}

For $n\ge4$ define
\[
  P(4) := p_8,\qquad
  P(n) := p_{P(n-1)} \quad (n\ge5),
\]
so that $P(n)$ is obtained from $8$ by $n-3$ nested applications of
the prime-index map:
\[
  P(n) = \underbrace{p_{p_{\cdots_{p_8}}}}_{n-3\ \text{times}}.
\]
For instance,
\[
\begin{aligned}
  P(4)&=p_8=19,\\
  P(5)&=p_{p_8}=p_{19}=67,\\
  P(6)&=p_{67}=331,\\
  P(7)&=p_{331}=2221.
\end{aligned}
\]

\begin{proposition}[Gutman--Ivi\'c, maximal case~\cite{gutman}]\label{prop:Matula-max}
For the Matula height function $H(p_i)=H(i)+1$, the largest number at
height $n$ is $P(n)$ for all $n\ge4$.  More precisely, for every
$m\in\mathcal{H}_n$ we have $m\le P(n)$, with equality if and only if
$m=P(n)$.
\end{proposition}

The proof relies on explicit bounds for the $k$th prime $p_k$.  We use
the following version of the Rosser--Schoenfeld inequalities.

\begin{lemma}[Rosser--Schoenfeld~\cite{RosserSchoenfeld}]\label{lem:RS}
There exists an absolute constant $k_0$ such that for all $k\ge k_0$,
\[
  k\bigl(\log k + \log\log k - 1\bigr)
  \;\le\; p_k
  \;\le\;
  k\bigl(\log k + \log\log k\bigr).
\]
\end{lemma}

We will not reproduce the proof of Lemma~\ref{lem:RS}; see
Rosser and Schoenfeld~\cite{RosserSchoenfeld}.  For small $k<k_0$ the
bounds can be checked directly.

\begin{proof}
We prove the result by induction.  For the proof only, set
\[
  P(0)=1,\qquad P(1)=2,\qquad P(2)=4,\qquad P(3)=8,
\]
in addition to the displayed definition of $P(n)$ for $n\ge4$.  Direct
inspection gives the finite base values
\[
\begin{array}{c|ccccc}
n&0&1&2&3&4\\
\hline
\max\mathcal H_n&1&2&4&8&19.
\end{array}
\]
The next few values are
\[
  P(5)=67,\qquad P(6)=331,\qquad P(7)=2221,\qquad P(8)=19577.
\]
These finite values absorb the small-index exceptions in the prime estimates
used below.

A direct finite computation of the height layers gives
\[
\max\mathcal H_5=67,\qquad
\max\mathcal H_6=331,\qquad
\max\mathcal H_7=2221,\qquad
\max\mathcal H_8=19577.
\]
These complete the base cases for the induction.

We now assume $n\ge9$ and that the theorem is known at all smaller
heights.

We first consider prime elements.  If $q\in\mathcal H_n$ is prime, then $q=p_i$ and
\[
  H(i)=H(p_i)-1=n-1.
\]
By the induction hypothesis, $i\le P(n-1)$, and hence
\[
  q=p_i\le p_{P(n-1)}=P(n).
\]
Thus, no prime element at height $n$ exceeds $P(n)$; equality in the prime case
forces $i=P(n-1)$.

We next consider composite elements.  Let $m\in\mathcal H_n$ be composite.  Choose a prime divisor $p$ of $m$ of
height $\ell\le n/2$; such a divisor exists because if every prime divisor had
height $>n/2$, then two prime divisors would already have total height exceeding
$n$.  Write $m=pu$, so $u\in\mathcal H_{n-\ell}$.  By induction,
\[
  p\le P(\ell),
  \qquad
  u\le P(n-\ell),
\]
and therefore
\begin{equation*}
  m\le P(\ell)P(n-\ell).
\end{equation*}
It remains to show
\begin{equation}\label{eq:Matula-supermult-short}
  P(\ell)P(n-\ell)\le P(n)
  \qquad
  \left(1\le \ell\le \left\lfloor\frac n2\right\rfloor,\ n\ge9\right).
\end{equation}

It remains to prove the required supermultiplicative estimate.  We use the standard bounds
\[
  p_k\ge k\log k,
  \qquad
  p_k\le k(\log k+\log\log k)
\]
in the ranges in which they are invoked; the finitely many smaller cases have
already been included in the base range above.

For $\ell\le3$, we have $P(\ell)\le8$.  Since $P(n-1)\ge P(8)=19577$ when
$n\ge9$, we have $\log P(n-1)>8$, and hence
\[
  P(n)=p_{P(n-1)}
  \ge P(n-1)\log P(n-1)
  >8P(n-1)
  \ge P(\ell)P(n-\ell).
\]
Thus, Equation~\eqref{eq:Matula-supermult-short} holds for $\ell\le3$.

Now assume $4\le \ell\le n/2$ and put $b=n-\ell$.  The lower bound gives
\[
  P(n)=p_{P(n-1)}\ge P(n-1)\log P(n-1),
\]
while the upper bound gives
\[
  P(\ell)=p_{P(\ell-1)}
  \le
  P(\ell-1)\bigl(\log P(\ell-1)+\log\log P(\ell-1)\bigr).
\]
By the induction hypothesis, $P(b)P(\ell-1)$ has height
\[
  b+(\ell-1)=n-1,
\]
and therefore
\begin{equation}\label{eq:Pk-1-dominates-product}
  P(n-1)\ge P(b)P(\ell-1).
\end{equation}
Also $n-1>\ell$, and the sequence $P(t)$ is strictly increasing since
$P(t+1)=p_{P(t)}>P(t)$.  Hence, $P(n-1)\ge P(\ell)$.  Applying the lower
bound to $P(\ell)=p_{P(\ell-1)}$ gives
\[
  P(n-1)\ge P(\ell)
  \ge P(\ell-1)\log P(\ell-1).
\]
Taking logarithms yields
\[
  \log P(n-1)
  \ge
  \log P(\ell-1)+\log\log P(\ell-1).
\]
Combining this with Equation~\eqref{eq:Pk-1-dominates-product} proves
Equation~\eqref{eq:Matula-supermult-short}.

Finally, we establish attainment and uniqueness.  The element $P(n)=p_{P(n-1)}$ itself has height $n$, so the maximum is attained.
The composite case above shows that every composite element of height $n$ is at
most $P(n)$, and since $P(n)$ is prime it cannot be equal to a composite
element.  In the prime case equality forces $i=P(n-1)$, hence the unique
maximal element is $P(n)$.

This completes the induction.
\end{proof}

Gutman and Ivi\'c also determined the smallest Matula number at height $n$.
Their original argument is phrased in terms of rooted trees; we give an
additive proof using only complete additivity of $H$ and a soft lower bound for
the $k$th prime.

\begin{lemma}[Prime lower bound]\label{lem:pk-gt-2k}
One has $p_k>2k$ for all $k\ge 6$.
\end{lemma}

\begin{proof}
By Dusart~\cite{dusart2018}, for all $k\ge 6$, 
\[
  p_k \ge k(\log k+\log\log k-1).
\]
Since $\log k+\log\log k-1>2$ for $k\ge 20$, it follows that $p_k>2k$; the remaining cases $6\le k\le 20$ are checked directly.
\end{proof}

\begin{proposition}[Gutman--Ivi\'c, minimal case~\cite{gutman}]\label{prop:Matula-min}
For the Matula height function $H(p_i)=H(i)+1$, the smallest number in
$\mathcal H_n$ is given, for $n>1$, by
\[
  m_{\min}(n) \;=\;
  \begin{cases}
    3\cdot 5^{(n-2)/3}, & n \equiv 2 \pmod 3,\\[4pt]
    5^{n/3},            & n \equiv 0 \pmod 3,\\[4pt]
    9\cdot 5^{(n-4)/3}, & n \equiv 1 \pmod 3.
  \end{cases}
\]
Moreover, this element is unique in $\mathcal H_n$.
\end{proposition}

\begin{proof}
For $n\ge 0$ define
\[
  a_n \;:=\; \min\{m\in\mathbb N:\ H(m)=n\}.
\]
A direct check gives
\[
  a_0=1,\quad a_1=2,\quad a_2=3,\quad a_3=5,\quad a_4=9.
\]
For $n>1$ set
\[
f(n)=
\begin{cases}
3\cdot 5^{(n-2)/3}, & n\equiv 2\pmod 3,\\[4pt]
5^{n/3},            & n\equiv 0\pmod 3,\\[4pt]
9\cdot 5^{(n-4)/3}, & n\equiv 1\pmod 3.
\end{cases}
\]
Then $a_n=f(n)$ holds for $2\le n\le 4$.

We first derive a recurrence.  For $n\ge 2$, complete additivity gives
\begin{equation}\label{eq:recurrence}
a_n=\min\Bigl\{\,p_{a_{n-1}},\ \min_{1\le k\le n-1} a_k a_{n-k}\Bigr\}.
\end{equation}
Indeed, $p_{a_{n-1}}$ has height $n$ because $H(a_{n-1})=n-1$, and
$a_ka_{n-k}$ has height $n$ for every $1\le k\le n-1$.  Thus, the right-hand
side gives admissible elements at height $n$.

Conversely, if $m\in\mathcal H_n$ is prime, then $m=p_i$ with
$H(i)=n-1$, so $i\ge a_{n-1}$ and $m\ge p_{a_{n-1}}$.  If $m$ is composite,
write $m=uv$ with $u,v>1$ and put $k=H(u)$.  Then $1\le k\le n-1$ and
$H(v)=n-k$, and so
\[
  m=uv\ge a_ka_{n-k}.
\]
This proves the recurrence.

We now proceed by induction.  Fix $n\ge 5$ and assume $a_m=f(m)$ for all $2\le m<n$.

We first show that the prime option is too large.  Since $n\ge 5$, we have $a_{n-1}\ge a_4=9\ge 6$, so by Lemma~\ref{lem:pk-gt-2k},
\[
  p_{a_{n-1}} \;>\; 2a_{n-1}.
\]
On the other hand, the composite integer $2a_{n-1}$ has height
$H(2a_{n-1})=H(2)+H(a_{n-1})=1+(n-1)=n$.

Hence, at height $n$
\[
  p_{a_{n-1}} \;>\;2a_{n-1}\,
\]
so the prime term in Equation~\eqref{eq:recurrence} cannot be minimal. Therefore,
\begin{equation}\label{eq:composite-only}
  a_n=\min_{1\le k\le n-1} a_k a_{n-k}
      =\min_{1\le k\le n-1} f(k)f(n-k),
\end{equation}
where we interpret $f(1):=a_1=2$.

We next evaluate the minimum.  Write $n=3b+r$ with $r\in\{0,1,2\}$ and $b\ge 1$.

\vskip 5pt\noindent {\tt Case 1:} $r=0$ ($n=3b$).
Write $k=3u+s$ with $s\in\{0,1,2\}$.  If $s=0$, then
\[
  f(k)f(n-k)=5^u\cdot 5^{b-u}=5^b=f(n).
\]
If $s=1$, then for $u=0$ one obtains
\[
  f(1)f(n-1)=2\cdot(3\cdot5^{b-1})=6\cdot5^{b-1}>5^b,
\]
while for $u\ge1$ a direct substitution gives a value larger than $5^b$.
If $s=2$, the boundary case $k=2$ gives
\[
  f(2)f(n-2)=3\cdot(9\cdot5^{b-2})=27\cdot5^{b-2}>5^b
\]
for $b\ge2$, and the remaining cases are again larger by the same substitution.
Thus, the minimum equals $5^b=f(n)$.

\vskip 5pt\noindent {\tt Case 2:} $r=2$ ($n=3b+2$).
Again write $k=3u+s$. If $s=0$ or $s=2$, then
\[
  f(k)f(n-k)=5^u\cdot (3\cdot 5^{b-u})=3\cdot 5^b=f(n)
\]
or
\[
  f(k)f(n-k)=(3\cdot 5^{u})\cdot 5^{b-u}=3\cdot 5^b=f(n),
\]
respectively.
If $s=1$ and $u=0$, then $f(1)f(n-1)=2\cdot (9\cdot 5^{b-1})=18\cdot 5^{b-1}>\!3\cdot 5^b$.
If $s=1$ and $u\ge 1$, then necessarily $n-k\neq 1$, and one computes
$f(k)f(n-k)\ge (9\cdot 5^{u-1})(9\cdot 5^{b-u-1})=81\cdot 5^{b-2}>\!3\cdot 5^b$ (when $b\ge 2$),
while the remaining small case $b=1$ is $n=5$, already covered by the base check.
Thus, the minimum equals $3\cdot 5^b=f(n)$.

\vskip 5pt\noindent {\tt Case 3:} $r=1$ ($n=3b+1$).
Write $k=3u+s$.  If $k\neq 1$ and $n-k\neq 1$, then a direct substitution from the
definition of $f$ gives in all three residue classes $s=0,1,2$ the identity
\[
  f(k)f(n-k)=9\cdot 5^{b-1}=f(n).
\]
If $k=1$ (or $n-k=1$), then $f(1)f(n-1)=2\cdot 5^b=10\cdot 5^{b-1}>\!9\cdot 5^{b-1}$.
Hence, the minimum equals $9\cdot 5^{b-1}=f(n)$.

\smallskip
Combining the three cases with Equation~\eqref{eq:composite-only} yields $a_n=f(n)$ for all $n\ge 5$,
completing the induction.

It remains to prove uniqueness.  The argument above shows that every minimizer in Equation~\eqref{eq:composite-only}
has value $f(n)$.  Although different splittings $k+(n-k)$ may attain the
minimum, the resulting product is always the same integer $f(n)$.  Hence, the
smallest element of $\mathcal H_n$ is unique and is equal to $f(n)$.
This proves uniqueness of $m_{\min}(n)$ in $\mathcal H_n$.
\end{proof}

\begin{remark}
Gutman and Ivi\'c used these extremal bounds, together with the tree interpretation of Matula numbers, to study the growth of $H(n)$ as $n\to\infty$.  They note that their results are obtained by graph-theoretic arguments and that ``it is not clear how Theorem~6 could be deduced by number-theoretical arguments, using only the defining properties (i) and (ii) of $f(n)$''~\cite[p.~141]{gutman}.  Propositions~\ref{prop:Matula-max} and~\ref{prop:Matula-min} provide such a number-theoretic derivation of the extremal bounds using only the recursion $H(p_i)=H(i)+1$, complete additivity, and standard estimates for $p_n$, without explicit reference to the underlying rooted trees.
\end{remark}

In light of this example, it is natural to ask what structure is represented by
$H(p_i)=H(i^2)+1$, or more generally by $H(p_i)=H(i^r)+1$.  For the square rule,
a prime at height $2j+1$ corresponds to an element of height $j$, and there are
no primes at even heights.  Thus, the primes at odd heights are in bijection with
the elements at half the preceding height, and the generating function is
\[
  \sum_{n=0}^{\infty}N_nq^n
  =
  \prod_{j\ge0}\frac{1}{(1-q^{2j+1})^{N_j}}.
\]
This corresponds to OEIS~A115593~\cite{OEIS}, the number of forests of rooted trees with
total weight $n$, where a node at height $k$ has weight $2^k$ and the root is
considered to have height $0$.  It is also a useful warning for the analytic
theory: recursive height rules can naturally produce profiles supported on a
proper arithmetic progression, so the minor-arc/aperiodicity hypothesis in
Section~\ref{sec6} is not automatic.

These examples suggest that recursive definitions of completely additive height
functions can produce combinatorial structures whose prime-height profiles are
not available in closed form.  In such cases the profile must often be computed
level by level, and the resulting weighted partition structure is part of the
object of study rather than an input.

\begin{example} \label{eg2.2}
Let $H(2)=1$ and $f(i) = (p_i - 1)/2$ for $i>1$. Then
\[
  H(p_i) = H\bigl((p_i - 1)/2\bigr) + 1 = H(p_i -1),
\]
which incidentally can also be written as $H(p_i)=H(\varphi(p_i))$ for $i>1$, where $\varphi$ is Euler's totient function. This sequence is described in OEIS~A006645~\cite{OEIS}.
In 1929, Pillai~\cite{pillai} introduced the height $H(n)$ of $n$ as the smallest positive integer $i$ such that the $i$th iterate of Euler's totient function at $n$ is $1$. This corresponds to the function $H(p_i) = H(p_i - 1) + 1$ but without $H$ being completely additive. It was shown in~\cite{shapiro1} that this $H$ was almost completely additive, failing to be completely additive only because $H(2p)=H(p)$. A modified version of $H$, which \emph{is} completely additive, was studied in~\cite{erdos} and is described here. The modified definition ensures $H(2p) = H(p) + 1$.
The difference in structure is simple: each odd number is shifted by one height from the original structure studied by Shapiro. We will see that this shift is common to an entire class of completely additive functions that emerge from the iterations of a multiplicative function.  The first few height layers are shown in Table~\ref{classes-shapiro}.\par
\begin{table}[ht]
\centering
\begin{tabular}[h]{| l | l | }
\hline\hline
$n$ & $\{m\mid H(m)=n\}$\cr
\hline
$5$ & $32, 34, 40, \mathbf{41}, 44, 46, \mathbf{47}, 48, 50, 51, 52, \mathbf{53}, 55, 56, 58, \mathbf{59}, 60, \mathbf{61}, 62, 65, 66,$\\
    & $ \mathbf{67}, 69, 70, \mathbf{71}, 72, \mathbf{73}, 74, 75, 76, 77, 78, \mathbf{79}, 84, 86, 87, 90, 91, 93, 95, 98, 99, 105, $\\
    & $ 108, \mathbf{109}, 111, 114, 117, 126, \mathbf{127}, 129, 133, 135, 147, 162, \mathbf{163}, 171, 189, 243 $\cr
$4$ & $ 16, \mathbf{17}, 20, 22, \mathbf{23}, 24, 25, 26, 28, \mathbf{29}, 30, \mathbf{31}, 33, 35, 36, \mathbf{37}, 38, 39, 42,$ \cr
    & $ \mathbf{43}, 45, 49, 54, 57, 63, 81$ \cr
$3$ & $ 8, 10, \mathbf{11}, 12, \mathbf{13}, 14, 15, 18, \mathbf{19}, 21, 27 $ \cr
$2$ & $ 4,\mathbf{5}, 6, \mathbf{7}, 9 $ \cr
$1$ & $ \mathbf{2}, \mathbf{3} $ \cr
$0$ & $1$ \cr
\hline
\end{tabular}
\caption{Numbers with height $\le 5$ for the modified $\varphi$-height. Primes are in bold.}\label{classes-shapiro}
\end{table} 
\end{example}

It is easy to show inductively that the smallest number at each height $n$ is
$2^n$ and the largest number is $3^n$.  If $\pi_n$ denotes the number of primes
at height $n$ and $N_n$ the total number of integers at height $n$, then
\[
  \sum_{n=0}^{\infty}N_nq^n
  =
  \prod_{n=1}^{\infty}\frac{1}{(1-q^n)^{\pi_n}}.
\]

Shapiro conjectured that there is at least one prime at each height, and that
remains an open question.  The computations in Section~\ref{sec8} suggest
growth roughly of the form $\pi_n\approx cB^n/n$ and $N_n\approx C B^n$, with
$B$ near $2.3$.  The fine structure of these Shapiro classes, and the close
connection between the height function and the prime-counting and prime-listing
functions $\pi(n)$ and $p_n$, was investigated in~\cite{balbhatnagar}.

These simple iterative definitions conceal a lot of structure. A natural question is how changes in the defining recursion alter the global height structure. If we consider the completely additive function defined by $H(p_i) = H(p_i - 1) +1$ (so without the special treatment of $2p$), we obtain the sequence set out in OEIS~A064097~\cite{OEIS}.

As noted in the comment there, if $d$ is the largest proper divisor of $n$, then $n/d = p$ where $p$ is a prime. Hence,
\[
\begin{aligned}
  H(dn/d)
  &=H(n/d)+H(d)\\
  &=H(p)+H(d)\\
  &=H(p-1)+1+H(d)\\
  &=H(d(p-1))+1\\
  &=H(n-d)+1.
\end{aligned}
\]
This proves that the function can be described as the number of iterations required for a number $n$ to reach $1$, such that at each iteration $n$ goes to $n-d$, where $d$ is the largest proper divisor of $n$.\par

Both Shapiro~\cite{shapiro2} and White~\cite{white} have examined the question of realizing height functions as iterations of multiplicative arithmetic functions. This case further suggests the connections with iterated functions leading to completely or almost completely additive functions, but we do not pursue this aspect in detail here.\par

Before we examine the number-theoretic aspects of these additive functions, including the distribution at each height, we look at one final example from the literature which is closely related to the previous example, but yields a rather different combinatorial structure.\par

\begin{example} \label{eg2.3}
Let $H(2)=1$ and $f(i) = (p_i + 1)/2$ for $i>1$. Then
\[
  H(p_i) = H\bigl((p_i + 1)/2\bigr) + 1 = H(p_i + 1),
\]
which can also be written as $H(p_i)= H(\sigma(p_i))$ for $i>1$, where $\sigma$ is the sum-of-divisors function. Dedekind studied the completely multiplicative function given by $f(p^\alpha)=p^{\alpha-1}(p+1)$, and Colin Defant~\cite{defant} studied a modification of this function which differs only at $p=2$, where it is given by $f(2^\alpha) = 2^{\alpha-1}$.\par
The original Dedekind function is a simple example of a completely multiplicative function whose iterates do not converge to $1$. Much like the Euler $\varphi$ function, the iterates of the modified function give rise to an almost completely additive function. The function studied here is the one modified further to yield a completely additive function.  Table~\ref{classes-dedekind} shows the first few height layers.\par
\begin{table}[ht]
\centering
\begin{tabular}[h]{| l | l | }
\hline\hline
$n$ & $\{m\mid H(m)=n\}$\cr
\hline
$6$ & $ 25, 27, \mathbf{29}, 30, 33, 34, 35, 36, \mathbf{37}, 38, 39, 40, \mathbf{41}, 42, \mathbf{43}, 44, 46, \mathbf{47}, 48, 49, $\\
    & $ 52, 56, \mathbf{61}, 62, 64 $\cr
$5$ & $ 15, \mathbf{17}, 18, \mathbf{19}, 20, 21, 22, \mathbf{23}, 24, 26, 28, \mathbf{31}, 32 $\cr
$4$ & $ 9, 10, \mathbf{11}, 12, \mathbf{13}, 14, 16$ \cr
$3$ & $ \mathbf{5}, 6, \mathbf{7}, 8 $ \cr
$2$ & $ \mathbf{3}, 4 $ \cr
$1$ & $ \mathbf{2}$ \cr
$0$ & $1$ \cr
\hline
\end{tabular}
\caption{Numbers with height $\le 6$ for the Dedekind–Defant type height. Primes are in bold.}\label{classes-dedekind}
\end{table} 
For this example one checks, by the same elementary induction used above, that
the smallest number at each height $n>2$ is given by
\[
  \begin{cases}
    5^{n/3},              & n \equiv 0\pmod 3,\\[4pt]
    9\cdot 5^{(n-4)/3},   & n \equiv 1\pmod 3,\\[4pt]
    3\cdot 5^{(n-2)/3},   & n \equiv 2\pmod 3,
  \end{cases}
\]
while the largest number is $2^n$.
\end{example}       

Many related examples can be obtained by altering the congruence rule in the
prime recursion.  For instance, one may define $H(2)=1$ and, for $i>1$, set
\[
 f(i) = \frac{p_i - 1}{3}  \quad\text{if }p_i \equiv 1\pmod 3,
 \qquad
 f(i) = \frac{p_i + 1}{3}  \quad\text{if }p_i \equiv 2\pmod 3.
\]
Analogous variants may be formed using congruence classes modulo $4$, $6$, or
other moduli.

\section{Squarefree Height Functions on Radical Classes}
\label{sec4}

We now pass to a squarefree variant of the preceding theory.  In
Sections~\ref{sec2} and~\ref{sec3}, a height function was an ordinary
completely additive function on \(\mathbb N\): prime powers contributed with
multiplicity, and the associated Euler product was the unrestricted
multipartition product
\[
  \prod_{k\ge1}(1-q^k)^{-\pi_k}.
\]
In the present section repeated prime factors are suppressed from the outset.
Thus, the relevant objects are radical classes, or equivalently their squarefree
representatives.  The counting problem changes from unrestricted
multipartitions to distinct multipartitions.

Define an equivalence relation on \(\mathbb N\) by
\[
  n\sim m
  \quad\text{if and only if}\quad
  \operatorname{rad}(n)=\operatorname{rad}(m),
\]
where \(\operatorname{rad}(n)\) is the product of the distinct prime divisors of
\(n\).  We shall always represent a class by its squarefree representative.

A squarefree height function is a function
\[
  \widehat H:\mathbb N\to\mathbb N_0
\]
which is constant on radical classes and additive over distinct prime factors:
\[
  \widehat H(1)=0,
  \qquad
  \widehat H(n)=\sum_{p\mid n}\widehat H(p).
\]
Equivalently,
\[
  \widehat H(n)=\widehat H(\operatorname{rad}(n)).
\]
Thus,
\[
  \widehat H(p^\alpha)=\widehat H(p)
  \qquad(\alpha\ge1),
\]
so a prime contributes at most once.

As before, we require
\[
  \widehat H(p)\ge1
\]
for every prime \(p\) in the structure, and we require finite prime fibers:
\[
  \#\{p:\widehat H(p)=k\}<\infty
  \qquad(k\ge1).
\]
Let
\[
  \widehat\pi_k=\#\{p:\widehat H(p)=k\},
\]
and let
\[
  \widehat N_n
  =
  \#\{r\text{ squarefree}:\widehat H(r)=n\}.
\]
Then
\[
  \sum_{n\ge0}\widehat N_nq^n
  =
  \prod_{k\ge1}(1+q^k)^{\widehat\pi_k}.
\]
Indeed, each prime of height \(k\) may occur either zero times or once in a
squarefree representative, and therefore contributes a factor \(1+q^k\).

For non-recursive assignments of prime heights, this gives the usual
distinct-part analog of the ordinary theory.  For recursive height
structures, however, one must be more careful.  In this section the recursion
is imposed directly on squarefree representatives.  We do not first construct
an ordinary completely additive height and then radicalize it.  Instead,
recursive prime-generation rules are applied only when the relevant input is
itself a squarefree representative built from the primes already present.  This
is the convention used throughout the rest of this section.

\begin{example}\label{eg:sf-distinct}
Let
\[
  \widehat H(p_i)=i.
\]
Then there is exactly one prime at each height.  The function is constant on
radical classes, since
\[
  \widehat H\Bigl(\prod_{j=1}^r p_{i_j}^{\alpha_j}\Bigr)
  =
  \sum_{j=1}^r i_j
  \qquad(\alpha_j>0).
\]
Therefore,
\[
  \sum_{n\ge0}\widehat N_nq^n
  =
  \prod_{n\ge1}(1+q^n).
\]
Thus, \(\widehat N_n\) is the number of partitions of \(n\) into distinct parts.  Table~\ref{tab:sf-distinct} shows the first few classes.
\begin{table}[ht]
\centering
\begin{tabular}{|l|l|}
\hline\hline
\(n\) & \(\{\, [m]\mid \widehat H(m)=n \,\}\)\\
\hline
\(6\) & \(\mathbf{[13]}, [21], [22], [30]\) \\
\(5\) & \(\mathbf{[11]}, [14], [15]\) \\
\(4\) & \(\mathbf{[7]}, [10]\) \\
\(3\) & \(\mathbf{[5]}, [6]\) \\
\(2\) & \(\mathbf{[3]}\) \\
\(1\) & \(\mathbf{[2]}\) \\
\(0\) & \([1]\) \\
\hline
\end{tabular}
\caption{Equivalence classes with height \(\le6\) for
\(\widehat H(p_i)=i\).  Primes are in bold.}
\label{tab:sf-distinct}
\end{table}
\end{example}

\begin{example}\label{eg:sf-phi}
Consider now the squarefree version of the Shapiro--Erd\H{o}s
\(\varphi\)-type recursion.  We begin with
\[
  \widehat H(2)=1.
\]
For odd primes \(p\), the squarefree analog of the relation
\(H(p)=H(p-1)\) is imposed as follows: a new prime \(p\) may enter the
structure only when
\[
  p-1
\]
is itself a squarefree product of primes already present.  In that case we set
\[
  \widehat H(p)=\widehat H(p-1).
\]
The squarefreeness condition is part of the recursion.  Thus, for instance,
\(5\) does not enter, because \(5-1=4\) is not squarefree.
A finite computation shows that the only primes that ever appear in this
squarefree structure are
\[
  2,\quad 3,\quad 7,\quad 43.
\]
This gives a finite certificate rather than an open-ended search.  Once the
candidate set is
\[
  S=\{2,3,7,43\},
\]
any further candidate prime must have the form
\[
  1+\prod_{p\in T}p
  \qquad(T\subseteq S),
\]
where the empty product gives the already known prime \(2\).  There are only
\(2^{|S|}=16\) such products to check, and none produces a new prime outside
\(S\).  Hence, no further primes occur.
The resulting squarefree representatives are exactly the squarefree products of
the primes in \(S\), with heights
\[
  \widehat H(2)=1,\qquad
  \widehat H(3)=1,\qquad
  \widehat H(7)=2,\qquad
  \widehat H(43)=4.
\]
Table~\ref{tab:sf-phi} lists the resulting squarefree classes.
\begin{table}[ht]
\centering
\begin{tabular}{|l|l|}
\hline\hline
\(n\) & \(\{\,m\text{ squarefree}:\widehat H(m)=n\,\}\)\\
\hline
\(6\) & \(258, 301\) \\
\(5\) & \(86, 129\) \\
\(4\) & \(42, \mathbf{43}\) \\
\(3\) & \(14, 21\) \\
\(2\) & \(6, \mathbf{7}\) \\
\(1\) & \(\mathbf{2}, \mathbf{3}\) \\
\(0\) & \(1\) \\
\hline
\end{tabular}
\caption{Squarefree representatives with height \(\le6\) in the squarefree
\(\varphi\)-height structure.  Primes are in bold.}
\label{tab:sf-phi}
\end{table}
\end{example}

\begin{example}\label{eg:sf-dedekind}
We finally consider the squarefree version of the Dedekind-type recursion from
Example~\ref{eg2.3}.  We take
\[
  \widehat H(2)=\widehat H(3)=1.
\]
The guiding relation is
\[
  p+1=\prod_{q\in T}q,
\]
where the product is squarefree and all primes \(q\) on the right have already
occurred.  For a newly generated prime \(p\), its height is the height of this
squarefree product:
\[
  \widehat H(p)=\widehat H(p+1).
\]
Thus, \(p+1\) must itself be a squarefree product of earlier primes.  This
condition is important: primes such as \(7\), for which \(7+1=8\), do not enter
the squarefree Dedekind-type structure, because \(8\) is not squarefree.
The first relations are
\[
  5+1=2\cdot3,\qquad
  29+1=2\cdot3\cdot5,\qquad
  173+1=2\cdot3\cdot29,\ \ldots.
\]
Equivalently,
\[
  2\cdot3-1=5,\qquad
  2\cdot3\cdot5-1=29,\qquad
  2\cdot3\cdot29-1=173,\ \ldots.
\]
Table~\ref{tab:sf-dedekind} lists the resulting squarefree classes.
\begin{table}[ht]
\centering
\begin{tabular}{|l|l|}
\hline\hline
\(n\) & \(\{\,m\text{ squarefree}:\widehat H(m)=n\,\}\)\\
\hline
\(6\) & \(145, \mathbf{173}, 174\) \\
\(5\) & \(58, 87\) \\
\(4\) & \(\mathbf{29}, 30\) \\
\(3\) & \(10, 15\) \\
\(2\) & \(\mathbf{5}, 6\) \\
\(1\) & \(\mathbf{2}, \mathbf{3}\) \\
\(0\) & \(1\) \\
\hline
\end{tabular}
\caption{Squarefree representatives with height \(\le6\) in the squarefree
Dedekind-type structure.  Primes are in bold.}
\label{tab:sf-dedekind}
\end{table}
It is natural to ask whether the primes generated by this squarefree
Dedekind-type rule are finite in number, as in the squarefree
\(\varphi\)-example, or whether the process continues indefinitely.  The
combinatorial setup allows us to list primes of this form iteratively.  The
first ten primes in the structure are
\[
  2,\ 3,\ 5,\ 29,\ 173,\ 5189,\ 300961,\ 4514429,\ 5386181,\ 161585429.
\]
Their heights, in the convention above, include
\[
  \widehat H(5189)=10,\qquad
  \widehat H(300961)=15,\qquad
  \widehat H(4514429)=18,
\]
\[
  \widehat H(5386181)=18,\qquad
  \widehat H(161585429)=22.
\]
Up to \(10^3\) there are \(5\) such primes, up to \(10^5\) there are \(6\), up
to \(10^6\) there are \(7\), and up to \(10^9\) there are \(11\).  This sparse
but continuing production of primes raises the question of whether infinitely
many primes occur in this squarefree structure.
To get a sense of the large-height behavior, consider the computed data for heights from \(96\) to \(100\) listed in Table~\ref{tab:sf-dedekind-large}.  Here \(\widehat N(n)\) denotes the
total number of squarefree representatives at height \(n\), and
\(\widehat\pi_n\) denotes the number of primes at that height.  The smallest
prime at height \(91\) is
\[
  859445547898845285802803723399409.
\]
\begin{table}[ht]
\centering
\[
\begin{array}{|c||c|c|c|c|}
\hline
n & \widehat N(n) & \widehat\pi_n
  & \widehat N(n)/\widehat N(n-1)
  & n\widehat\pi_n/\widehat N(n) \\
\hline
96 & 133461 & 1751 & 1.133 & 1.260 \\ \hline
97 & 151302 & 1901 & 1.134 & 1.219 \\ \hline
98 & 171567 & 2154 & 1.134 & 1.230 \\ \hline
99 & 194421 & 2308 & 1.133 & 1.176 \\ \hline
100 & 220134 & 2742 & 1.132 & 1.246 \\ \hline
\end{array}
\]
\caption{\(\widehat N(n)\) and \(\widehat\pi_n\) for \(n=96,\dots,100\) in the
squarefree Dedekind-type structure.}
\label{tab:sf-dedekind-large}
\end{table}
These data suggest that this squarefree structure may approach a stable
asymptotic regime: the ratios
\[
  \widehat N(n)/\widehat N(n-1)
\]
and the scaled densities
\[
  n\widehat\pi_n/\widehat N(n)
\]
appear to stabilize over the computed range.  We record this only as
computational evidence.  It motivates the broader question of whether such
squarefree recursive height structures admit a finite-height probabilistic
theory, but no asymptotic theorem is asserted here.
\end{example}

\section{The Average Order of the Height Functions}\label{sec5}

For a height function $H$ we study the summatory function
\[
  F(x) \;:=\; \sum_{n\le x} H(n),
\]
which, for completely additive $H$, admits a clean reduction to the values of
$H$ on the primes.

\subsection{Reduction to the Prime Values}

\begin{proposition}[Prime reduction]\label{prop:prime-reduction}
Let $H$ be completely additive. Then for every $x\ge 2$,
\begin{equation}\label{eq:Fx-prime-decomp}
  F(x)
  \;=\;
  \sum_{p^\alpha\le x} H(p)\,\Bigl\lfloor \frac{x}{p^\alpha}\Bigr\rfloor
  \;=\;
  x\sum_{p\le x} H(p)\!\!\sum_{\substack{\alpha\ge1\\ p^\alpha\le x}}\frac{1}{p^\alpha}
  \;+\;
  O\!\Bigl(\sum_{p^\alpha\le x} H(p)\Bigr).
\end{equation}
Consequently,
\begin{equation}\label{eq:Fx-main-term}
  F(x)
  =
  x\sum_{p\le x}\frac{H(p)}{p-1}
  +
  O\!\left(
    \sum_{p\le x} H(p)\frac{\log x}{\log p}
  \right)
  +
  O\!\left(\sum_{p\le x}H(p)\right).
\end{equation}
In particular,
\[
  F(x)
  =
  x\sum_{p\le x}\frac{H(p)}{p-1}
  +
  O\!\left(
    \log x\sum_{p\le x}H(p)
  \right).
\]
\end{proposition}

\begin{proof}
Writing $n=\prod_p p^{v_p(n)}$, complete additivity gives
\[
H(n)=\sum_p v_p(n)H(p)=\sum_p\sum_{\alpha=1}^{v_p(n)}H(p)=\sum_{p^\alpha\mid n}H(p).
\]
Summing over $n\le x$ and exchanging summations yields
\[
F(x)=\sum_{n\le x}\sum_{p^\alpha\mid n}H(p)
=\sum_{p^\alpha\le x}H(p)\,\#\{n\le x:\,p^\alpha\mid n\}
=\sum_{p^\alpha\le x}H(p)\Bigl\lfloor \frac{x}{p^\alpha}\Bigr\rfloor,
\]
which is the first identity in Equation~\eqref{eq:Fx-prime-decomp}.  Writing
$\lfloor x/p^\alpha\rfloor=x/p^\alpha+O(1)$ gives the displayed form of
Equation~\eqref{eq:Fx-prime-decomp}.  The floor-error term is bounded by
\[
  \sum_{p^\alpha\le x}H(p)
  \le
  \sum_{p\le x}H(p)\left\lfloor\frac{\log x}{\log p}\right\rfloor.
\]
Also,
\[
  \sum_{\substack{\alpha\ge1\\p^\alpha\le x}}\frac1{p^\alpha}
  =
  \frac1{p-1}+O\!\left(\frac1x\right),
\]
so the truncation error contributes $O(\sum_{p\le x}H(p))$.  This yields
Equation~\eqref{eq:Fx-main-term}.
\end{proof}

\begin{remark}\label{rem:pminus1}
The replacement of $(p-1)^{-1}$ by $p^{-1}$ must be handled with some care.
In general
\[
  \frac1{p-1}=\frac1p+O\!\left(\frac1{p^2}\right),
\]
so the resulting change in the main sum is bounded by
\[
  x\sum_{p\le x}\frac{H(p)}{p^2}.
\]
This is lower order in the examples where we make the replacement, but it is not
uniformly $O(x)$ without additional hypotheses on $H(p)$.
\end{remark}

\subsection{Examples}

We briefly recall the average-order behavior for the principal examples in this paper and
compare it with the growth of the associated weighted-partition multiplicities $N_n$.

First, consider the sum-of-primes height $H(p)=p$.  Here complete additivity gives
\[
  H(n)=\sum_p v_p(n)p,
\]
the sum of prime factors counted with multiplicity.  Alladi and Erd\H{o}s showed
that
\[
  F(x)=\sum_{n\le x}H(n)\sim C_1\frac{x^2}{\log x}.
\]
This matches the heuristic from Equation~\eqref{eq:Fx-main-term}:
$x\sum_{p\le x}H(p)/p\sim x\pi(x)\sim x^2/\log x$.  The associated partition
problem is the problem of partitions into prime parts, since the profile has one
generator at height $p$ for each prime $p$; this is not one of the polynomial
profiles treated directly by Theorem~\ref{thm:inverse-growth} and has additional
slowly varying factors.

Second, consider polynomial heights in the primes, with $H(p)=p^k$.  If $H(p)=p^k$, then
\[
\sum_{p\le x}\frac{H(p)}p
=
\sum_{p\le x}p^{k-1}
\sim \frac{x^k}{k\log x}.
\]
Thus, the leading prime sum in Equation~\eqref{eq:Fx-main-term} has size
$x^{k+1}/\log x$, and one obtains the order of magnitude
\[
  F(x)\asymp \frac{x^{k+1}}{\log x}
\]
by the standard estimates for additive functions with nonnegative prime values.
On the partition side the allowed part sizes are $k$th powers of primes, not all
$k$th powers.

Third, consider Matula heights.  For the Matula height, the associated multiplicities $N_n$ (rooted trees) grow exponentially,
so $\log N_n\asymp n$ (Otter~\cite{otter}).  On the additive side,
Bret\`eche and Tenenbaum~\cite{breteche} showed that
\[
  F(x)\asymp x\log x.
\]
This lies outside the polynomial regime treated by the inverse-growth theorem of
Section~\ref{sec6}.

Finally, consider iterated $\varphi$ heights.  The difficulty of obtaining average-order information for the Shapiro height
is already visible at the level of primes.  Let $\phi_0(n)=n$ and
$\phi_{k+1}(n)=\phi(\phi_k(n))$, and write $k(n)$ for the least $k$ such that
$\phi_k(n)=1$. Shapiro observed that a slight parity adjustment of $k(n)$ yields a completely
additive function, which we denote here by $H_\varphi(n)$.  Erd\H{o}s, Granville, Pomerance, and Spiro make this explicit: $H_\varphi$ is completely
additive and, on primes, satisfies the recursion
\[
  H_\varphi(p)=H_\varphi(p-1)\qquad(p>2).
\]
Thus, the height of primes is governed by the arithmetic of the shifted set
$\{p-1:p\ \text{prime}\}$, and the prime-height profile
$\pi_k=\#\{p:H_\varphi(p)=k\}$ is entangled with fine distributional
information on primes in arithmetic progressions.
Indeed, the strongest known normal/average-order results for $H_\varphi(n)$ in
\cite{erdos} are proved conditionally on Elliott--Halberstam type hypotheses,
under which one obtains an asymptotic
$\frac1x\sum_{n\le x}H_\varphi(n)=\alpha\log x+o(\log x)$ together with
concentration estimates implying that $H_\varphi(n)$ has normal order
$\alpha\log n$.
This provides a concrete explanation for why Shapiro-type height questions lie
well beyond the ``designed profile'' setting treated elsewhere in this paper.

These examples motivate the inverse-growth principle developed in the next
section: in the polynomial-profile regime, the growth of the prime-height profile
determines the dominant analytic pole, and under the full hypotheses of
Meinardus' theorem this pole controls the stretched-exponential growth of the
multiplicities $N_n$.  

\section{Inverse Growth for Height Functions via Weighted Partitions}
\label{sec6}

The weighted-partition identity of Proposition~\ref{thm:height-multipartition}
places the height multiplicities in the classical setting of Euler products with
weighted parts.  In this section we use this theory in one direction: we show how
regular growth of the prime-height profile gives the pole data of a Dirichlet
series, and how, under the full hypotheses of Meinardus' theorem, this pole data
yields a stretched-exponential law for the coefficients.

It is important to separate these two steps.  A cumulative estimate for the
profile determines the dominant pole of the Dirichlet series, but it does not by
itself imply the minor-arc condition in Meinardus' theorem.  The latter condition
rules out lattice obstructions, such as profiles supported only on even heights.

Throughout this section, $H:\mathbb N\to\mathbb N_0$ denotes a completely
additive height function as in Section~\ref{sec2}, and
\[
  \pi_k:=\#\{p\text{ prime}:H(p)=k\},
  \qquad
  N_n:=\#\{m\ge1:H(m)=n\}.
\]
Thus,
\[
  \sum_{n\ge0}N_nq^n
  =
  \prod_{k\ge1}(1-q^k)^{-\pi_k}.
\]

\subsection{Weighted Partitions and a Theorem of Meinardus}

Let $(b_k)_{k\ge1}$ be a sequence of nonnegative real numbers, not identically
zero, and set
\[
  F(q)=\prod_{k\ge1}(1-q^k)^{-b_k}=\sum_{n\ge0}a_nq^n,
  \qquad
  B(s)=\sum_{k\ge1}\frac{b_k}{k^s},
\]
\[
  g(\tau)=\sum_{k\ge1}b_ke^{-k\tau}\qquad(\Re\tau>0).
\]
We write $s=\sigma+it$.  We use the following standard form of Meinardus'
theorem~\cite{meinardus}; the minor-arc hypothesis is included explicitly
because it is essential for the applications below.

\begin{theorem}[Meinardus~\cite{meinardus}]
\label{thm:Meinardus}
Suppose that for some $\alpha>0$ the following hypotheses hold.
\begin{enumerate}
  \item[\textup{(M1)}]
  $B(s)$ converges for $\sigma>\alpha$ and admits a meromorphic continuation to
  a half-plane $\sigma\ge -C_0$ for some $0<C_0<1$, where it is holomorphic
  except for a single simple pole at $s=\alpha$ with residue $A>0$.
  \item[\textup{(M2)}]
  There is a constant $c>0$ such that
  \[
    B(\sigma+it)=O\bigl((1+|t|)^c\bigr)
  \]
  uniformly for $\sigma\ge -C_0$ as $|t|\to\infty$.
  \item[\textup{(M3)}]
  \textup{(Minor-arc/aperiodicity condition.)} There exist constants
  $\varepsilon>0$ and $C_3>0$ such that, writing $\tau=y+2\pi ix$,
  \[
    \Re g(y+2\pi ix)-g(y)\le -C_3y^{-\varepsilon}
  \]
  uniformly for
  \[
    y^{1+\varepsilon/2}\le |x|\le \tfrac12
  \]
  as $y\to0^+$.
\end{enumerate}
Then the coefficients have the usual Meinardus asymptotic
\[
  a_n\sim C_1n^\kappa
  \exp\!\left(C_2n^{\frac{\alpha}{\alpha+1}}\right)
\]
for some $C_1>0$, where
\[
  \kappa=\frac{B(0)-1-\tfrac12\alpha}{\alpha+1}
\]
and
\begin{equation*}
  C_2=
  \left(1+\frac1\alpha\right)
  \bigl(A\Gamma(\alpha+1)\zeta(\alpha+1)\bigr)^{1/(\alpha+1)}.
\end{equation*}
In particular,
\[
  \log a_n\sim C_2n^{\alpha/(\alpha+1)}.
\]
\end{theorem}

\begin{remark}[Why the minor-arc condition is necessary]
\label{rem:M3-necessary}
Condition \textup{(M3)} is not implied by the pole and vertical-growth
conditions.  Let $b_k=1$ for $k$ even and $b_k=0$ for $k$ odd.  Then
\[
  F(q)=\prod_{m\ge1}(1-q^{2m})^{-1},
\]
so $a_n=0$ for every odd $n$.  However
\[
  B(s)=2^{-s}\zeta(s)
\]
has a simple pole at $s=1$ and satisfies the usual vertical-growth condition in
strips.  What fails is precisely the minor-arc condition: the support lies on a
proper lattice.  This is not a contrived pathology; recursive height rules can naturally produce profiles supported on a proper arithmetic progression, as in the square-rule variant \(H(p_i)=H(i^2)+1\) discussed above.
\end{remark}

\subsection{From Prime-Height Growth to Pole Data}

\begin{lemma}[Pole data from cumulative profile growth]
\label{lem:pole-from-partial-sums}
Let $(\pi_k)_{k\ge1}$ be nonnegative and suppose that
\begin{equation}\label{eq:pi-partial-sum}
  \Pi(x):=\sum_{k\le x}\pi_k
  =
  \frac{C}{\alpha}x^\alpha+O\bigl(x^{\alpha-\eta}\bigr)
\end{equation}
for some $C>0$, $\alpha>0$, and $\eta>0$.  Then
\[
  B_H(s):=\sum_{k\ge1}\frac{\pi_k}{k^s}
\]
converges for $\Re s>\alpha$ and continues meromorphically to
$\Re s>\alpha-\eta$, with a single simple pole in that region, at $s=\alpha$,
of residue $C$.  Moreover, for every $\delta\in(0,\eta)$ one has
$B_H(\sigma+it)=O(1+|t|)$ uniformly for
$\sigma\ge\alpha-\eta+\delta$.
\end{lemma}

\begin{proof}
For $\Re s>\alpha$, partial summation gives
\[
  B_H(s)=s\int_1^\infty x^{-s-1}\Pi(x)\,dx.
\]
Writing $\Pi(x)=\frac{C}{\alpha}x^\alpha+E(x)$ with
$E(x)=O(x^{\alpha-\eta})$, the main term contributes
\[
  \frac{Cs}{\alpha(s-\alpha)},
\]
which has residue $C$ at $s=\alpha$.  The integral
$s\int_1^\infty x^{-s-1}E(x)\,dx$ is holomorphic for
$\Re s>\alpha-\eta$ and satisfies the stated vertical bound on every closed
sub-half-plane $\sigma\ge\alpha-\eta+\delta$.
\end{proof}

\begin{remark}[What the cumulative estimate does and does not supply]
\label{rem:cumulative-not-full-meinardus}
The estimate in Equation~\eqref{eq:pi-partial-sum} supplies the dominant pole and residue of
$B_H(s)$ in the half-plane reached by the error term.  It does not automatically
supply the full half-plane continuation required in \textup{(M1)}, nor does it
supply the minor-arc condition \textup{(M3)}.  Thus, applications of
Meinardus' theorem require either a direct verification of \textup{(M1)}--
\textup{(M3)} or an independent reference guaranteeing those hypotheses.  In
particular, a power law for $\Pi(x)$ is compatible with periodic support, as in
Remark~\ref{rem:M3-necessary}.
\end{remark}

\subsection{An Inverse-Growth Theorem for Height Functions}

\begin{theorem}[Conditional inverse growth]
\label{thm:inverse-growth}
Let $H$ be a completely additive height function with profile $(\pi_k)$ and
height multiplicities $(N_n)$.  Suppose
\[
  \sum_{n\ge0}N_nq^n=\prod_{k\ge1}(1-q^k)^{-\pi_k}.
\]
Assume that the cumulative profile satisfies Equation~\eqref{eq:pi-partial-sum} and that
the associated weights $b_k=\pi_k$ satisfy the full hypotheses
\textup{(M1)}--\textup{(M3)} of Theorem~\ref{thm:Meinardus}, with the pole at
$s=\alpha$ having residue $A=C$.  Then there exist $C_1>0$ and
$\kappa\in\mathbb R$ such that
\[
  N_n\sim C_1n^\kappa
  \exp\!\left(C_2n^{\frac{\alpha}{\alpha+1}}\right),
\]
where
\begin{equation*}
  C_2=
  \left(1+\frac1\alpha\right)
  \bigl(C\Gamma(\alpha+1)\zeta(\alpha+1)\bigr)^{1/(\alpha+1)}.
\end{equation*}
In particular,
\[
  \log N_n\sim C_2n^{\alpha/(\alpha+1)}.
\]
\end{theorem}

\begin{proof}
The weighted-partition identity gives the Euler product with weights
$b_k=\pi_k$.  By hypothesis, these weights satisfy the full hypotheses of
Theorem~\ref{thm:Meinardus}, and the residue at the dominant pole is $C$.
The result is therefore exactly Theorem~\ref{thm:Meinardus} applied to this
Euler product.
\end{proof}

\begin{remark}
The adjective ``conditional'' is essential.  The profile estimate determines
the predicted exponent and constant once the analytic hypotheses hold, but it
cannot by itself rule out a lattice obstruction.  For example, if the profile is
supported on even heights, an all-$n$ asymptotic is impossible because every odd
coefficient vanishes.
\end{remark}

\subsection{From the Growth of $H(p_i)$ to the Profile Exponent}

The exponent $\alpha$ in Theorem~\ref{thm:inverse-growth} is the exponent in the
cumulative prime-height growth
\[
  \Pi(x)=\sum_{k\le x}\pi_k=\#\{p:H(p)\le x\}.
\]
If the prime heights are arranged so that $H(p_i)$ grows regularly with the
prime index $i$, then $\Pi(x)$ may be read off by inversion: heuristically, if
$H(p_i)\asymp i^d$, then $\Pi(x)\asymp x^{1/d}$, so $\alpha=1/d$ and the
stretched-exponential exponent becomes
\[
  \frac{\alpha}{\alpha+1}=\frac{1}{d+1}.
\]
In Section~\ref{sec7} we make this inversion canonical by fixing, for a given
profile, the sequential realization in which prime heights are nondecreasing in
prime order.

\subsection{Identification of the Parameters from $(N_n)$}

\begin{corollary}[Identification under Meinardus hypotheses]
\label{thm:inverse-converse}
Let $H$, $\pi_k$, and $N_n$ be as above.  Suppose that the Dirichlet series
\[
  B_H(s)=\sum_{k\ge1}\frac{\pi_k}{k^s}
\]
satisfies the full hypotheses \textup{(M1)}--\textup{(M3)} of
Theorem~\ref{thm:Meinardus}, with rightmost singularity a simple pole at
$s=r>0$ of residue $A>0$.  Suppose moreover that
\[
  N_n\sim C_1n^\kappa\exp(C_2n^\beta)
\]
for some $\beta\in(0,1)$, $C_1,C_2>0$, and $\kappa\in\mathbb R$.  Then
\[
  \beta=\frac{r}{r+1},
  \qquad
  r=\frac{\beta}{1-\beta}.
\]
Writing $\alpha=r$, one has
\[
  C_2=
  \left(1+\frac1\alpha\right)
  \bigl(A\Gamma(\alpha+1)\zeta(\alpha+1)\bigr)^{1/(\alpha+1)}.
\]
Moreover,
\[
  \sum_{k\le x}\pi_k\sim\frac{A}{\alpha}x^\alpha.
\]
\end{corollary}

\begin{proof}
Theorem~\ref{thm:Meinardus} gives
\[
  \log N_n\sim
  \left(1+\frac1r\right)
  \bigl(A\Gamma(r+1)\zeta(r+1)\bigr)^{1/(r+1)}n^{r/(r+1)}.
\]
Comparison with the assumed asymptotic forces $\beta=r/(r+1)$ and hence
$r=\beta/(1-\beta)$.  The displayed formula for $C_2$ follows by substituting
$r=\alpha$.  Finally, since $\pi_k\ge0$ and $B_H$ is holomorphic on
$\Re s\ge\alpha$ except for the simple pole at $s=\alpha$, the
Wiener--Ikehara theorem for Dirichlet series with nonnegative coefficients gives
$\sum_{k\le x}\pi_k\sim(A/\alpha)x^\alpha$; see, for example,
\cite[Theorem~III.1]{tenenbaum} or \cite[Ch.~V, \S4]{korevaar}.
\end{proof}

\begin{remark}
This is a conditional converse: it presupposes the full analytic hypotheses and
a single rightmost simple pole.  In the heuristic model $H(p_i)\asymp i^d$, it
recovers $\alpha=1/d$ and $\beta=1/(d+1)$.
\end{remark}

\subsection{Examples}

\begin{example}[Ordinary partitions]
\label{ex:ordinary-partitions-inverse}
Let $\pi_k\equiv1$.  Then $B_H(s)=\zeta(s)$, $\alpha=1$, and $A=1$.  The
classical minor-arc condition is satisfied in this case; equivalently, this is
the usual Hardy--Ramanujan partition problem, treated for instance in
Andrews~\cite[Ch.~5--6]{andrews}.  Theorem~\ref{thm:Meinardus} recovers
\[
  p(n)\sim \frac{1}{4n\sqrt3}\exp\!\left(\pi\sqrt{\frac{2n}{3}}\right).
\]
Here $\kappa=-1$ and
\[
  C_2=2\bigl(\Gamma(2)\zeta(2)\bigr)^{1/2}
  =\pi\sqrt{\frac23}.
\]
\end{example}

\begin{example}[Plane partitions]
\label{ex:plane-partitions-inverse}
Let $\pi_k=k$.  Then
\[
  B_H(s)=\sum_{k\ge1}k^{1-s}=\zeta(s-1),
\]
so $\alpha=2$ and $A=1$.  The full hypotheses are classical in this case; this
is the MacMahon plane-partition product, whose asymptotics go back to
Wright~\cite{wright}.  Theorem~\ref{thm:Meinardus} gives
\[
  \log N_n\sim \frac32\bigl(2\zeta(3)\bigr)^{1/3}n^{2/3},
\]
in agreement with Wright's plane-partition asymptotic.  This example is a useful
check on the constant: replacing $(1+1/\alpha)$ by $(\alpha+1)$ would make the
constant too large by a factor of $2$ when $\alpha=2$.
\end{example}

These examples illustrate the role of the full hypotheses.  The pole of the
Dirichlet series determines the exponent and constant once the minor-arc
condition is known, but a cumulative power law alone is not enough to guarantee
an all-$n$ coefficient asymptotic.

\section{Canonical Sequential Realizations and Average Order in the Polynomial Regime}\label{sec7}

\subsection{Why a Canonical Realization Is Needed in the Profile Regime}
In the polynomial (non-explosive) regime of \Cref{sec6}, the essential input for the \emph{layer sizes}
$N_n=\#\{m:\height(m)=n\}$ is the \emph{multipartition profile}
\[
\pi_k := \#\{p \ \text{prime} : \height(p)=k\},\qquad k\ge 1,
\]
since it governs the Euler product
\[
\sum_{n\ge 0}N_n q^n=\prod_{k\ge 1}(1-q^k)^{-\pi_k}.
\]
However, unlike recursive models (e.g.\ Matula-type heights) where the recursion fixes $\height$ uniquely,
a profile $\{\pi_k\}$ generally admits many completely additive realizations, obtained by permuting
which primes receive which heights.  In this section we select a canonical representative and
compute its average order using the prime-sum reduction of \Cref{sec5}.
This complements \Cref{sec6}, which concerns the asymptotic growth of the layer sizes $N_n$ as a
weighted partition problem; here we study the \emph{average order} of a canonical height realizing the
same profile.

\subsection{Definition of the Canonical Sequential Model}
Define cumulative counts
\[
\Pi_k := \sum_{j\le k}\pi_j.
\]

\begin{definition}[Canonical sequential realization]\label{def:Hseq}
Let $(\pi_k)_{k\ge1}$ be a profile with $\pi_k<\infty$ for every $k$ and
\[
  \sum_{k\ge1}\pi_k=\infty.
\]
Let
\[
  \Pi_k=\sum_{j\le k}\pi_j,\qquad \Pi_0=0,
\]
and let $p_1<p_2<\cdots$ be the primes.  Define
$\height_{\mathrm{seq}}$ on primes by blocks:
\[
\height_{\mathrm{seq}}(p_i)=k
\quad\text{whenever}\quad
\Pi_{k-1}< i\le \Pi_k,
\]
and extend to all $n=\prod p_i^{a_i}$ by complete additivity:
\[
\height_{\mathrm{seq}}(n)=\sum_i a_i\,\height_{\mathrm{seq}}(p_i).
\]
We call the resulting height function the \emph{canonical sequential realization}
of the profile.
\end{definition}

Equivalently, $\height_{\mathrm{seq}}$ is the unique realization of the given profile whose values
on primes are nondecreasing in the natural prime order.

\subsection{Prime Sums and Average Order}
Let
\[
F_{\mathrm{seq}}(x):=\sum_{n\le x}\height_{\mathrm{seq}}(n),
\qquad
S_{\mathrm{seq}}(x):=\sum_{p\le x}\frac{\height_{\mathrm{seq}}(p)}{p}.
\]
By Equation~\eqref{eq:Fx-main-term} in \Cref{sec5}, the average-order problem reduces
to estimating prime sums:
\[
F_{\mathrm{seq}}(x)
=
x\sum_{p\le x}\frac{\height_{\mathrm{seq}}(p)}{p-1}
+
O\!\left(
  \sum_{p\le x}\height_{\mathrm{seq}}(p)\frac{\log x}{\log p}
\right).
\]
Thus, \(S_{\mathrm{seq}}(x)\) determines the scale of the leading prime sum,
although the displayed error term is generally of the same order in the
polynomial-profile regime.

\subsection{Polynomial Cumulative Profiles: $\Pi_k\sim Ck^\alpha$}

\begin{lemma}[Inverting a polynomial profile]\label{lem:invertPi}
Assume $\Pi_k\sim Ck^\alpha$ as $k\to\infty$ with $C>0$ and $\alpha>0$, and that $\Pi_k$ is eventually
strictly increasing.  Let $k(i)$ be the unique integer such that $\Pi_{k(i)-1}< i\le \Pi_{k(i)}$.
Then
\[
k(i)\sim \left(\frac{i}{C}\right)^{1/\alpha}\qquad (i\to\infty),
\]
and hence $\height_{\mathrm{seq}}(p_i)=k(i)\sim (i/C)^{1/\alpha}$.
\end{lemma}

\begin{proof}
By definition, $\Pi_{k(i)-1}<i\le \Pi_{k(i)}$.  Since $\Pi_k\sim Ck^\alpha$ and $\Pi_k$ is eventually
increasing, the inverse function satisfies $k(i)\sim (i/C)^{1/\alpha}$ by standard monotone inversion:
for any $\varepsilon>0$ and large $k$,
\[
(C-\varepsilon)k^\alpha \le \Pi_k \le (C+\varepsilon)k^\alpha,
\]
which implies
\[
\left(\frac{i}{C+\varepsilon}\right)^{1/\alpha}\le k(i)\le
\left(\frac{i}{C-\varepsilon}\right)^{1/\alpha}
\]
for all large $i$, and letting $\varepsilon\to 0$ gives the claim.
\end{proof}

\begin{lemma}[A logarithmic-integral sum]\label{lem:logint-sum}
Fix $\beta>0$. Then, as $m\to\infty$,
\[
\sum_{2\le i\le m}\frac{i^{\beta-1}}{\log i}
\;=\;
\frac{m^\beta}{\beta\log m}\,(1+o(1)).
\]
\end{lemma}

\begin{proof}
Fix $\varepsilon\in(0,1)$.  The contribution from $2\le i\le m^{1-\varepsilon}$
is
\[
  O\!\left(\frac{m^{\beta(1-\varepsilon)}}{\log m}\right)
  =
  o\!\left(\frac{m^\beta}{\log m}\right).
\]
On the remaining range $m^{1-\varepsilon}<i\le m$ we have
$\log i=(1+O(\varepsilon))\log m$, uniformly.  Hence,
\[
\sum_{m^{1-\varepsilon}<i\le m}\frac{i^{\beta-1}}{\log i}
=
\frac{1+O(\varepsilon)}{\log m}
\sum_{m^{1-\varepsilon}<i\le m}i^{\beta-1}.
\]
Since
\[
  \sum_{i\le m}i^{\beta-1}\sim \frac{m^\beta}{\beta},
\]
and the omitted initial segment is negligible, we obtain
\[
  \sum_{2\le i\le m}\frac{i^{\beta-1}}{\log i}
  =
  (1+O(\varepsilon))
  \frac{m^\beta}{\beta\log m}.
\]
Letting $\varepsilon\to0$ proves the claim.
\end{proof}

\begin{theorem}[Prime-harmonic growth for the canonical sequential model]\label{thm:seq-prime-sum}
Assume $\Pi_k\sim Ck^\alpha$ with $C>0$ and $\alpha>0$, and let $m=\pi(x)$.
Then
\[
S_{\mathrm{seq}}(x)=\sum_{p\le x}\frac{\height_{\mathrm{seq}}(p)}{p}
\;=\;
\frac{\alpha}{C^{1/\alpha}}\cdot \frac{m^{1/\alpha}}{\log m}\,(1+o(1)).
\]
In particular, using $\pi(x)\sim x/\log x$,
\[
S_{\mathrm{seq}}(x)\asymp \frac{x^{1/\alpha}}{(\log x)^{1+1/\alpha}}.
\]
\end{theorem}

\begin{proof}
Write $m=\pi(x)$ and sum over prime indices:
\[
S_{\mathrm{seq}}(x)=\sum_{i\le m}\frac{\height_{\mathrm{seq}}(p_i)}{p_i}.
\]
By \Cref{lem:invertPi}, $\height_{\mathrm{seq}}(p_i)\sim (i/C)^{1/\alpha}$.
Also $p_i\sim i\log i$ (Prime Number Theorem in the form $p_i\sim i\log i$), hence
\[
\frac{\height_{\mathrm{seq}}(p_i)}{p_i}
\sim
\frac{C^{-1/\alpha} i^{1/\alpha}}{i\log i}
=
C^{-1/\alpha}\frac{i^{(1/\alpha)-1}}{\log i}.
\]
Let $\beta:=1/\alpha$. Summing and applying \Cref{lem:logint-sum} gives
\[
S_{\mathrm{seq}}(x)
\sim
C^{-1/\alpha}\sum_{i\le m}\frac{i^{\beta-1}}{\log i}
\sim
C^{-1/\alpha}\cdot \frac{m^\beta}{\beta\log m}
=
\frac{\alpha}{C^{1/\alpha}}\cdot\frac{m^{1/\alpha}}{\log m}.
\]
Finally, $\pi(x)\sim x/\log x$ yields the displayed scale in $x$.
\end{proof}

\begin{corollary}[Average-order scale from Equation~\eqref{eq:Fx-main-term}]
\label{cor:seq-avgorder-scale}
Assume $\Pi_k\sim Ck^\alpha$ as in \Cref{thm:seq-prime-sum}, and let
$F_{\mathrm{seq}}(x)=\sum_{n\le x}\height_{\mathrm{seq}}(n)$.  Then the leading
prime-sum term has scale
\[
  x\sum_{p\le x}\frac{\height_{\mathrm{seq}}(p)}{p}
  \asymp
  \frac{x^{1+1/\alpha}}{(\log x)^{1+1/\alpha}}.
\]
Moreover, the error term in Equation~\eqref{eq:Fx-main-term} has the same order of
magnitude.  Consequently
\[
  F_{\mathrm{seq}}(x)
  \asymp
  \frac{x^{1+1/\alpha}}{(\log x)^{1+1/\alpha}}.
\]
The reduction determines the order of magnitude in the polynomial-profile
regime, but not an asymptotic constant for $F_{\mathrm{seq}}(x)$.
\end{corollary}

\begin{proof}
Let $m=\pi(x)$.  By \Cref{thm:seq-prime-sum},
\[
  xS_{\mathrm{seq}}(x)
  \asymp
  x\frac{m^{1/\alpha}}{\log m}
  \asymp
  \frac{x^{1+1/\alpha}}{(\log x)^{1+1/\alpha}},
\]
using $m\sim x/\log x$.

For the error term in Equation~\eqref{eq:Fx-main-term}, use
$\height_{\mathrm{seq}}(p_i)\asymp i^{1/\alpha}$ and
$p_i\sim i\log i$.  Then
\[
\sum_{p\le x}\height_{\mathrm{seq}}(p)\frac{\log x}{\log p}
\asymp
\log x\sum_{i\le m}\frac{i^{1/\alpha}}{\log i}.
\]
Applying \Cref{lem:logint-sum} with $\beta=1+1/\alpha$ gives
\[
\sum_{i\le m}\frac{i^{1/\alpha}}{\log i}
\asymp
\frac{m^{1+1/\alpha}}{\log m}.
\]
Since $\log x\asymp\log m$, this error term is
\[
  \asymp m^{1+1/\alpha}
  \asymp
  \frac{x^{1+1/\alpha}}{(\log x)^{1+1/\alpha}}.
\]
For the lower bound, use the exact identity
\[
  F_{\mathrm{seq}}(x)
  =
  \sum_{p^\nu\le x}
  \height_{\mathrm{seq}}(p)
  \left\lfloor \frac{x}{p^\nu}\right\rfloor.
\]
Keeping only the terms with $\nu=1$ and $p\le x/2$ gives
\[
  F_{\mathrm{seq}}(x)
  \ge
  \sum_{p\le x/2}
  \height_{\mathrm{seq}}(p)
  \left\lfloor \frac{x}{p}\right\rfloor
  \ge
  \frac{x}{2}
  \sum_{p\le x/2}\frac{\height_{\mathrm{seq}}(p)}p,
\]
which has the same order of magnitude by \Cref{thm:seq-prime-sum}.

Finally, replacing $(p-1)^{-1}$ by $p^{-1}$ changes the leading prime sum by
\[
  x\sum_{p\le x}\frac{\height_{\mathrm{seq}}(p)}{p^2},
\]
which is lower order than $xS_{\mathrm{seq}}(x)$ by the same calculation as in
Remark~\ref{rem:pminus1}.  This proves the asserted order of magnitude and
explains why no asymptotic constant follows from the reduction alone.
\end{proof}

\subsection{Two Benchmark Profiles: Partitions and Plane Partitions}

\begin{corollary}[Ordinary partitions profile]\label{cor:seq-partitions}
For the ordinary partition profile $\pi_k\equiv 1$ (so $\Pi_k=k$, hence $\alpha=1$ and $C=1$),
the canonical sequential model satisfies
\[
S_{\mathrm{seq}}(x)\asymp \frac{x}{(\log x)^2},
\qquad
F_{\mathrm{seq}}(x)\asymp \frac{x^2}{(\log x)^2}.
\]
\end{corollary}

\begin{proof}
Here $\Pi_k=k$, so \Cref{thm:seq-prime-sum} gives
$S_{\mathrm{seq}}(x)\asymp \pi(x)/\log\pi(x)\asymp x/(\log x)^2$.
Then \Cref{cor:seq-avgorder-scale} yields $F_{\mathrm{seq}}(x)\asymp x^2/(\log x)^2$.
\end{proof}

\begin{corollary}[Plane partitions profile]\label{cor:seq-plane-partitions}
For the plane partition profile $\pi_k=k$ (so $\Pi_k=\tfrac12k(k+1)\sim \tfrac12 k^2$, hence
$\alpha=2$ and $C=\tfrac12$), the canonical sequential model satisfies
\[
S_{\mathrm{seq}}(x)\sim 2\sqrt{2}\,\frac{\sqrt{\pi(x)}}{\log \pi(x)}
\asymp \frac{x^{1/2}}{(\log x)^{3/2}},
\]
and
\[
F_{\mathrm{seq}}(x)\asymp x\,S_{\mathrm{seq}}(x)\asymp \frac{x^{3/2}}{(\log x)^{3/2}}.
\]
\end{corollary}

\begin{proof}
Apply \Cref{thm:seq-prime-sum} with $\alpha=2$ and $C=1/2$ to obtain the stated main term for
$S_{\mathrm{seq}}(x)$.  The scale for $F_{\mathrm{seq}}(x)$ then follows from
\Cref{cor:seq-avgorder-scale}.
\end{proof}

\subsection{Beyond Polynomial Profiles and Transition to Shapiro}
The polynomial profile assumption $\Pi_k\sim Ck^\alpha$ is precisely the regime treated cleanly by the
inverse-growth theory of \Cref{sec6}, and it is also the regime in which the canonical sequential model
admits a transparent prime-sum analysis as above.  Profiles supported sparsely on arithmetic sets (for
instance, $\pi_k$ supported on primes) still admit a canonical sequential realization, but sharp average
order in such cases typically requires finer control of large prime-factor statistics than we pursue here.
We now turn to the recursively defined Shapiro totient height, where prime structure is not imposed by a
profile but emerges from the recursion itself.

\section{Computational Observations for the Shapiro Height Structure}
\label{sec8}

Up to this point our analysis of height functions has been essentially
\emph{vertical}: we have either averaged $H(n)$ as $n\le x$ grows, or we have
encoded the distribution of prime heights into weighted partition problems.
In this section we study the structure \emph{horizontally}, at fixed height.
Concretely, we fix a height $h$ and examine the $h$th layer of the Shapiro
height structure: the multiplicities
\[
  N_h := \#\{n\ge1 : H(n)=h\},\qquad
  P_h := \#\{p\text{ prime} : H(p)=h\},
\]
and the distribution of the primes at height $h$ by size.

We focus on the modified Shapiro height function from Example~\ref{eg2.2}.
Throughout this section $H$ denotes this completely additive height.  The
material in this section is empirical.  The tables and figures are included to
record stable numerical patterns and to formulate conjectures; they are not used
in the proofs of the main theorems, and no new theorem about $\pi(x)$ is claimed.

\subsection{Exponential Growth of $N_h$ and $P_h$}

The first striking feature of the data is the regular growth of $N_h$ and
$P_h$ with $h$.  Table~\ref{table:NhPhData} lists $N_h$ and $P_h$ for
$1\le h\le17$, together with the ratios $N_h/N_{h-1}$ and the scaled ratios
\[
  \frac{P_h\cdot h}{P_{h-1}\cdot(h-1)}.
\]

\begin{table}[ht]
\centering
\[
\begin{array}{| c | c| c || c| c| c |}
\hline
h & N_h & N_h/N_{h-1} & h & P_h & \dfrac{P_h\cdot h}{P_{h-1}\cdot(h-1)} \\
\hline
1  & 2       & 2.00 & 1  & 2       & -    \cr \hline
2  & 5       & 2.50 & 2  & 2       & 2.00 \cr \hline
3  & 11      & 2.20 & 3  & 3       & 2.25 \cr \hline
4  & 26      & 2.36 & 4  & 6       & 2.67 \cr \hline
5  & 59      & 2.27 & 5  & 12      & 2.50 \cr \hline
6  & 137     & 2.32 & 6  & 23      & 2.30 \cr \hline
7  & 312     & 2.28 & 7  & 46      & 2.33 \cr \hline
8  & 719     & 2.30 & 8  & 94      & 2.34 \cr \hline
9  & 1651    & 2.30 & 9  & 198     & 2.37 \cr \hline
10 & 3816    & 2.31 & 10 & 424     & 2.38 \cr \hline
11 & 8757    & 2.29 & 11 & 854     & 2.22 \cr \hline
12 & 20202   & 2.31 & 12 & 1859    & 2.37 \cr \hline
13 & 46440   & 2.30 & 13 & 3884    & 2.26 \cr \hline
14 & 106957  & 2.30 & 14 & 8362    & 2.32 \cr \hline
15 & 245989  & 2.30 & 15 & 17837   & 2.29 \cr \hline
16 & 566561  & 2.30 & 16 & 38977   & 2.33 \cr \hline
17 & 1303968 & 2.30 & 17 & 84188   & 2.30 \cr \hline
\end{array}
\]
\caption{$N_h$ and $P_h$ for $h=1,\dots,17$ for the modified Shapiro height.}
\label{table:NhPhData}
\end{table}

The ratios $N_h/N_{h-1}$ stabilize rapidly near $2.30$, suggesting an
exponential law $N_h\approx C_N B^h$ with $B\approx 2.3$.  Likewise,
the quantities $hP_h/(P_{h-1}(h-1))$ stabilize near a constant, in line with
a heuristic of the form $P_h\approx C_P B^h/h$.
Motivated by this, we record the following empirical regularity.

\begin{conjecture}[Exponential growth in height]\label{conj:Nh-Ph}
For the modified Shapiro height function, there exist constants $B>1$,
$C_N>0$, and $C_P>0$ such that
\[
  N_h \sim C_N B^h,
  \qquad
  P_h \sim C_P \frac{B^h}{h}
\]
as $h\to\infty$.
Numerically, the data up to $h=17$ are well-fitted by $B\approx 2.3$.
\end{conjecture}

\subsection{Lognormal Fits at Fixed Height}

\begin{figure}[ht]
  \centering
  \includegraphics[width=\textwidth]{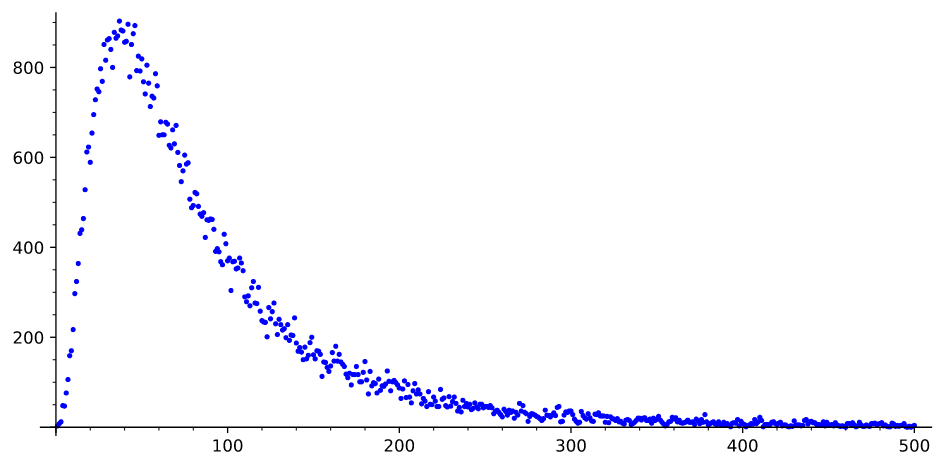}
  \caption{Distribution of primes at height $17$ for the modified Shapiro height function.
  Each point represents the number of primes in an interval of length $25000$ in
  the prime variable $p$ after shifting the origin to $2^{17}$; the overlaid curve corresponds to a normal fit in the variable $\log p$.}
  \label{fig:lognormal-primes}
\end{figure}

A second pattern appears when we examine the distribution of primes at a fixed
height.  For each $h$ we collect all primes with $H(p)=h$ and study the
distribution of their logarithms.  Figure~\ref{fig:lognormal-primes} shows the
height-$17$ data as a histogram in the prime variable $p$, using bins of width
$25000$ after shifting the origin to $2^{17}$; the overlaid curve corresponds to
a Gaussian fit for $\log p$.

To quantify this behavior, we computed for each height $h$ the mean and
standard deviation of $\log p$ over all primes with $H(p)=h$, and then
normalized these by $h$ and $\sqrt{h}$ respectively.  The resulting statistics
are shown in Table~\ref{table:height-stats}.

\begin{table}[ht]
\centering
\[
\begin{array}{| c || c | c |}
\hline
h & \text{Mean}(\log p)/h & \text{SD}(\log p)/h^{0.5} \\
\hline
1  & 0.8900 & 0.2000 \cr \hline
2  & 0.8888 & 0.1319 \cr \hline
3  & 0.8879 & 0.1419 \cr \hline
4  & 0.8392 & 0.1535 \cr \hline
5  & 0.8580 & 0.1761 \cr \hline
6  & 0.8487 & 0.1805 \cr \hline
7  & 0.8453 & 0.1726 \cr \hline
8  & 0.8464 & 0.1766 \cr \hline
9  & 0.8444 & 0.1656 \cr \hline
10 & 0.8504 & 0.1805 \cr \hline
11 & 0.8487 & 0.1773 \cr \hline
12 & 0.8487 & 0.1764 \cr \hline
13 & 0.8486 & 0.1771 \cr \hline
14 & 0.8486 & 0.1779 \cr \hline
15 & 0.8486 & 0.1763 \cr \hline
16 & 0.8486 & 0.1774 \cr \hline
17 & 0.8486 & 0.1772 \cr \hline
\end{array}
\]
\caption{Normalized mean and standard deviation of $\log p$ for primes with
$H(p)=h$, for $1\le h\le17$.}
\label{table:height-stats}
\end{table}

From $h\approx 10$ onwards, the normalized quantities are essentially
constant:
\[
  \frac{\mathbb{E}(\log p \mid H(p)=h)}{h}
  \approx 0.8486,\qquad
  \frac{\mathrm{SD}(\log p \mid H(p)=h)}{\sqrt{h}}
  \approx 0.1771.
\]
This motivates the following ``height-wise central limit'' conjecture.

\begin{conjecture}[Height-wise central limit law for primes]
\label{conj:lognormal}
There exist constants $\mu,\sigma>0$ such that, for every fixed real $t$,
\[
  \frac{1}{P_h}
  \#\left\{
    p\text{ prime}: H(p)=h,\ 
    \frac{\log p-\mu h}{\sigma\sqrt h}\le t
  \right\}
  \longrightarrow
  \Phi(t)
\]
as $h\to\infty$, where $\Phi$ is the standard normal distribution function.
Equivalently, the primes at height $h$ are conjecturally distributed like a
lognormal family with parameters $(\mu h,\sigma^2h)$ in the variable $p$.
Numerically, the data up to $h=17$ are consistent with
$\mu\approx0.8486$ and $\sigma\approx0.1771$.
\end{conjecture}

We emphasize that Conjecture~\ref{conj:lognormal} is supported only by finite
computations; at present we have no theoretical mechanism for proving it.

\subsection{A Height-Stratified Proxy for \texorpdfstring{$\pi(x)$}{pi(x)}}

Assuming Conjectures~\ref{conj:Nh-Ph} and~\ref{conj:lognormal}, one can build a
simple height-stratified proxy for the prime-counting function $\pi(x)$.
Fix $x>0$ and consider the contribution from a single height $h$.  The total
number of primes at height $h$ is heuristically $P_h\approx C_PB^h/h$, and by
Conjecture~\ref{conj:lognormal} the logarithms of these primes are
approximately $N(\mu h,\sigma^2 h)$.  In the displayed proxy below, the fitted
multiplicative constant $C_P$ is absorbed into the normalization; on the
computed range it is close to $1$.  Summing the resulting contributions over
$1\le h\le \lceil \log x/\log 2\rceil$ leads to
\begin{equation}\label{primefit}
  \hat{\pi}(x)
  \;=\;
  \sum_{h=1}^{\left\lceil{\frac{\log x}{\log 2}}\right\rceil}
    \frac{B^{h}}{h}
    \left(
      \frac{1}{2}
      + \frac{1}{2}\,\mathrm{erf}\!\left(
        \frac{\log x - \mu h}{\sqrt{2h}\,\sigma}
      \right)
    \right),
\end{equation}
where in practice we take the fitted values
\[
  B\approx 2.3,\qquad \mu\approx 0.8486,\qquad \sigma\approx 0.1771.
\]

Combining Conjectures~\ref{conj:Nh-Ph} and~\ref{conj:lognormal} gives the
height-stratified proxy $\hat{\pi}(x)$ in Equation~\eqref{primefit}, calibrated from
finite layer data.  As a coarse consistency check, comparison with the tabulated
value of $\pi(10^{29})$ gives:
\[
\begin{aligned}
\pi(10^{29})&=1.520698109714\times 10^{27},\\
\hat{\pi}(10^{29})&=1.519637333970\times 10^{27},\\
\frac{10^{29}}{\log(10^{29})}&=1.497567178976\times 10^{27}.
\end{aligned}
\]
We do not regard Equation~\eqref{primefit} as a proposed approximation method for
$\pi(x)$.  Its role is more limited: it shows that the two empirical layer laws
recorded above are numerically compatible with the global scale of prime
counting far outside the fitted range.

Thus, this section supplies concrete conjectures, rather than theorem-level
results, linking three observed phenomena: (i) exponential growth of $N_h$ and
$P_h$, (ii) Gaussian-like behavior of $\log p$ at fixed height, and (iii) a
height-stratified model whose aggregate scale is consistent with prime counts.

\section{Extensions, Universality, and Outlook}
\label{sec9}

The main structural mechanism of this paper is the passage
\[
  H \longmapsto \bigl(\pi_k\bigr)_{k\ge1},
  \qquad
  \pi_k=\#\{p:\ H(p)=k\},
\]
followed by the weighted-partition Euler product
\[
  \sum_{n\ge0} N_n q^n \;=\; \prod_{k\ge1}(1-q^k)^{-\pi_k},
  \qquad
  N_n=\#\{m:\ H(m)=n\}.
\]
Thus, once the prime-height profile is known, the height multiplicities are
determined by a purely combinatorial partition model.  This viewpoint also
extends verbatim to any multiplicative system with unique factorization
(e.g.\ ideals in $\mathcal{O}_K$, or $\mathbb{F}_q[x]$), by replacing primes by
atoms; we do not pursue these extensions here, since the central questions in
this paper already appear over $\mathbb{N}$.

From an analytic perspective, the key input is the growth of the cumulative
profile
\[
  \Pi(x)\;:=\;\sum_{k\le x}\pi_k.
\]
This leads naturally to two broad regimes.

The first case, which we call Regime~I, is the polynomial-profile case.  If $\Pi(x)$ has power-law growth (as in the benchmark multipartition profiles
$\pi_k\equiv1$ and $\pi_k=k$), then this growth supplies the pole and residue of
the profile Dirichlet series
\[
  B_H(s)=\sum_{k\ge1}\frac{\pi_k}{k^s}.
\]
When the remaining hypotheses of Meinardus' theorem, including the minor-arc
condition, are satisfied, Theorem~\ref{thm:inverse-growth} gives a precise
\emph{inverse-growth} law for $\log N_n$ with an explicit constant.  Thus, the
polynomial profile determines the predicted exponent, while the analytic
hypotheses rule out lattice obstructions such as support on a proper arithmetic
progression.  Section~\ref{sec7} shows that, in this same polynomial regime, one
can also attach to a given profile a canonical sequential realization and derive
corresponding average-order consequences via the additive-function reduction of
Section~\ref{sec5}.

The second case, which we call Regime~II, is the exponential-profile case.  If $\Pi(x)$ grows essentially exponentially (equivalently, $\pi_k$ grows like
$B^k$ up to polynomial factors), then $B_H(s)$ typically diverges for all $s$
and the Dirichlet--Meinardus approach breaks down.  In this regime the Euler
product usually has finite radius of convergence, and one expects the dominant
singularity of $\prod_{k\ge1}(1-q^k)^{-\pi_k}$ to govern the leading exponential
growth of $N_n$.  Developing a general theory for such profiles remains an
interesting direction for future work.

The examples treated here organize naturally along this divide.  Ordinary and
plane partitions sit in Regime~I and serve as benchmark profiles for the
inverse-growth theorem.  The Matula height, driven by the intrinsic recursion
$H_M(p_i)=H_M(i)+1$, belongs to an exponential class, reflecting the exponential
growth of rooted-tree counts.  The Shapiro totient height also appears
exponential in the tested range (Section~\ref{sec8}), but with strong arithmetic
constraints coming from the factorization of $p-1$; it therefore provides a
natural test case for how ``prime dynamics'' can shape behavior beyond the
universal polynomial picture.

A structural limitation of the partition correspondence is that the map
$H\mapsto(\pi_k)$ forgets substantial information.  From the point of view of
$\prod_{k\ge1}(1-q^k)^{-\pi_k}$, only the \emph{counts} of primes at each height
matter; congruence information, correlations, and finer arithmetic structure are
invisible.  This is a strength in Regime~I, where coarse asymptotics depend only
on $(\pi_k)$, but it also highlights where genuinely arithmetic phenomena must
enter: refined statistics live in the fibers over a fixed profile and cannot be
recovered from $(\pi_k)$ alone.

Two directions suggested most directly by the present work are the following.
First, one would like a systematic ``Regime~II'' analog of Section~\ref{sec6},
based on singularity analysis of Euler products with exponentially growing
weights.  Second, the Shapiro computations in Section~\ref{sec8} indicate that
height layers can carry meaningful \emph{horizontal} structure in addition to
the vertical multiplicities $N_h$: in our data the prime sizes at fixed height
exhibit an approximately Gaussian profile in $\log p$ after rescaling.  This
suggests studying conditional statistics in the microcanonical ensemble
conditioned on $H=n$, and comparing them to the canonical (Boltzmann-weighted)
models naturally associated to the Euler product.

\section*{Acknowledgements}
I would like to thank Gaurav Bhatnagar for many helpful discussions and for insights that shaped some of the work in this paper.

\end{document}